\documentclass[hyperref]{ctexart}
\usepackage{geometry} 
\usepackage{helvet}
\usepackage{amsmath, amsfonts, amssymb}
\usepackage[english]{babel}

\usepackage{url}      
\usepackage{bm}      
\usepackage{multirow}
\usepackage{booktabs}
\usepackage{cite}
\usepackage{lipsum}
\usepackage{graphicx}
\usepackage{graphics}
\usepackage{subfigure}
\usepackage{epstopdf}
\usepackage{float}
\usepackage{indentfirst}
\usepackage{algorithm}
\usepackage{algorithmicx}  
\usepackage{algpseudocode} 

\usepackage{fancyhdr} 
\usepackage{hyperref} 
\hypersetup{colorlinks, bookmarks, unicode} 

\numberwithin{equation}{section}

\newtheorem{Theorem}{Theorem}[section]
\newtheorem{Lemma}{Lemma}[section]
\newtheorem{Conjecture}{Conjecture}[section]
\newtheorem{Property}{Property}[section]
\newtheorem{Proof}{Proof}[section]
\newtheorem{Remark}{Remark}[section]

\title{\textbf{Riemann problem for constant flow with single-point heating source}}
\author{\sffamily Changsheng Yu, \sffamily Chengliang Feng, \sffamily Zhiqiang Zeng, \sffamily Tiegang Liu}

\begin{document}
	\maketitle
	
\begin{abstract}
	This work focuses on the Riemann problem of Euler equations with global constant initial conditions and a single-point heating source, which comes from the physical problem of heating one-dimensional inviscid compressible constant flow. In order to deal with the source of Dirac delta-function, we propose an analytical frame of double classic Riemann problems(CRPs) coupling, which treats the fluids on both sides of the heating point as two separate Riemann problems and then couples them. Under the double CRPs frame, the solution is self-similar, and only three types of solution are found. The theoretical analysis is also supported by the numerical simulation. Furthermore, the uniqueness of the Riemann solution is established with some restrictions on the Mach number of the initial condition.
\end{abstract}
	


\section{Introduction}
\begin{figure}[h]
	\label{figure device}
	\centering
	\includegraphics[height=3.5cm]{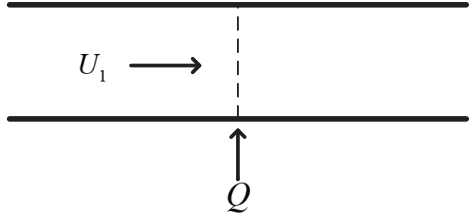}\par
	\caption{one-dimensional constant flow with heat addition from a single point}
\end{figure}

The Riemann problem of the one-dimensional inviscid compressible flow with global constant initial conditions and a singe-point heating source is studied in this paper. The heating point is located at $x=0$. The governing equations is given by 
\begin{equation}
	\label{governing equations}
	\frac{\partial U}{\partial t}+\frac{\partial F}{\partial x}=S,
\end{equation}
where
\[
U=
\begin{pmatrix}
	\rho\\\rho u\\E
\end{pmatrix},
F=
\begin{pmatrix}
	\rho u\\\rho u^2+p\\(E+p)u
\end{pmatrix},
S=
\begin{pmatrix}
	0\\0\\Q \delta(x)
\end{pmatrix}.
\]
Here, $\rho$, $p$ and $E$ denote the thermodynamical variables: density, pressure and total energy, respectively. $u$ is velocity. $Q>0$ is the heat flux per unit time added to the flow. $\delta(x)$ is the Dirac delta-function. The source means that $Q$ heat is added to the flow from the heating point per unit time in the physical sense. We assume that the fluid is polytropic ideal and the equaiotn of state is given by
\[ p=(\gamma-1)\rho e,\quad 1\leq \gamma\leq 3.  \]
where $\gamma$ is the ratio of specific heats and $e$ is the internal energy. The initial condition is
\begin{equation}
	\label{initial condition}
	U(x,0)\equiv U_1 =(\rho_1, \rho_1u_1, E_1),
\end{equation}
where $\rho_1$, $u_1$ and $E_1$ are constant. There is no loss of generality in assuming $u_1>0$. The subscript "1" means the initial state in this paper. The physical problem described by the Riemann problem (\ref{governing equations}) and (\ref{initial condition}) is that $Q$ heat is added to the one-dimensional constant flow per unit time, as shown in Figure\ref{figure device}.

The subscripts "-" and "+" denote the limiting states upstream and downstream of the heating point,  respectively. Write
\[
U_-(t)=U(0-,t), \quad U_+(t)=U(0+,t).
\]

The source term implies a jump of the energy flux across the heating point.

\begin{subequations}
	\label{balance equation}
	\begin{align}
		\rho_-(t)u_-(t)&=\rho_+(t)u_+(t)\\
		\rho_-(t){u_-(t)}^2+p_-(t)&=\rho_+(t){u_+(t)}^2+p_+(t)\\
		(E_-(t)+p_-(t))u_-(t)+Q&=(E_+(t)+p_+(t))u_+(t).
	\end{align}
\end{subequations}

If the velocity at the heating point is zero, thermal convection does not take effect, then the heat addition has no effect on the flow. Thus in the remainder of this paper we assume that
\[
u_-(t)\ne 0, \quad u_+(t)\ne 0.
\]

We only deal with the situation of constant $U_-(t)$ and $U_+(t)$. Omitting the time parameter $t$, we can get the following equation from (\ref{balance equation}).
\[
\left(\frac{1}{2}{u_-}^2+h_- \right)(1+k) =\frac{1}{2}{u_+}^2+h_+,
\]
where $h$ is enthalpy. $k=\frac{Q}{\rho_-(t)u_-(t)\left(\frac{1}{2}{u_-}^2+h_- \right)}$ is called the heating parameter and is assumed to be a constant parameter.

The heating equations are defined by
\begin{equation}
	\label{balance equations2}
	\left\lbrace 
	\begin{aligned}
		&\rho_-u_-=\rho_+u_+\\
		&\rho_-{u_-}^2+p_-=\rho_+{u_+}^2+p_+\\
		&\left(\frac{1}{2}{u_-}^2+h_- \right)(1+k) =\frac{1}{2}{u_+}^2+h_+
	\end{aligned}.
	\right. 
\end{equation}
The solution of the heating equations (\ref{balance equations2}) corresponds to a steady solution of the equations (\ref{governing equations}).

The mathematical model of Figure\ref{figure device} can be applied to the study of one-dimensional condensation problem (see \cite{2007Condensation,2002On}). The condensation of the vapor leads to the release of latent heat and therefore has heating effects on the carrier fluid. Previous researches on condensation problems have revealed some information about the solution of Riemann problem (\ref{governing equations}) and (\ref{initial condition}). The solution of the heating equations (\ref{balance equations2}) has already been available in \cite{2007Condensation,delale1993the,2002On,schnerr2005unsteadiness}. In \cite{schnerr2005unsteadiness} Schnerr made a good summary of the properties of the solution to (\ref{balance equations2}) at the case of $M_->1$, which is typical for condensation problems. Schnerr showed that there are two solutions of the heating equations (\ref{balance equations2}), one called the shock solution, which reduces to identity when $Q=0$, and the other called the weak solution, which reduces to the adiabatic normal shock solution when $Q=0$. For the heat addition of subsonic flow, Schnerr predicted the appearance of the unsteady solution, but did not do a more detailed analysis. For the unsteady solution of this Riemann problem, Dongen et al.\cite{2002On} studied the unsteady effects of the heat addition, and proposed three possible solution structures. Those three structures are determined by the Mach numbers of the fluid around the heating point. However, a complete and rigorous theoretical proof is lacking at present. On the other hand, numerical simulation as an effective tool has been employed to explore the wave patterns of the heating problem. Chengwan et al.\cite{cheng2010on} used the ASCE method\cite{luo2006on} to numerically verify the above three structures by the simulation of the onset of condensation in a slender Laval nozzle. In the wet nitrogen condensation problem of the Laval nozzle, Chengwan showed the transition between the three structures by adjusting the humidity. To our knowledge, there is no effective theoretical method to study the exact solution of Riemann problem (\ref{governing equations}) and (\ref{initial condition}).

Conservative hyperbolic equations with source terms can be transformed into non-conservative hyperbolic equations without source terms. Adopting the following form for the $\delta$-function
\[
\delta(x)=\frac{\partial H(x)}{\partial x}, \quad 
H(x)=\begin{cases}
	0, \quad x<0\\
	1, \quad x> 0
\end{cases},
\]
where $H(x)$ is called Heaviside function, one can rewrite the Riemann problem (\ref{governing equations}) and (\ref{initial condition}) into following non-conservative form.
\begin{equation}
	\label{augmented system}
	\frac{\partial \widetilde{U}}{\partial t}+\widetilde{A}\frac{\partial \widetilde{U}}{\partial x}=0,
\end{equation}
where
\[
\widetilde{U}=
\begin{pmatrix}
	\rho\\\rho u\\E\\h
\end{pmatrix},
\widetilde{A}=
\begin{pmatrix}
	0&1&0&0\\
	\frac{1}{2}(\gamma-3)u^2&(3-\gamma)u&\gamma-1&0\\
	\frac{1}{2}(\gamma-2)u^3-\frac{c^2u}{\gamma-1}&\frac{3-2\gamma}{2}u^2+\frac{c^2}{\gamma-1}&\gamma u&-Q\\
	0&0&0&0
\end{pmatrix},
\]
where $c$ is the speed of sound. The initial condition for the augmented equations is 
\[\widetilde{U}(x,0)=\begin{pmatrix}
	\rho_0&\rho_0 u_0&E_0&H(x)
\end{pmatrix}^T.\]

By defining proper entropy solution for non-conservative equations (\ref{augmented system}), one can analyzed the solutions of the original Riemann problem and design appropriate numerical methods to make numerical simulation (see \cite{abgrall2010comment,gosse2001well,greenberg1997analysis}). Applications of that method include the fluid in a nozzle with discontinuous cross-sectional area(see \cite{kroner2005numerical,lefloch2003the,thanh2009the}) and the shallow water equations with discontinuous topography(see \cite{ALCRUDO2001643,bernetti2008exact,lefloch2007the,pares2019the,thanh2013numerical}). The augmented equations (\ref{augmented system}) give rise to an additional linearly degenerated characteristic field and an additional stationary discontinuity. Although the augmented equations do not contain source terms, the solution process of the Riemann problem is very complex and often problem-related due to the lack of conservation or strict hyperbolicity.

The Riemann problem of homogeneous Euler equations with piecewise constant initial conditions is called classical Riemann problem(CRP). Altough the addition of the source term does not lead to the loss of self-similarity, the singularity of the source term has a significant effect on the structure of the solution.
The Riemann solution of the Euler equations with smooth source terms, which are often called generalized Riemann problem (GRP, see \cite{benartzi1984a,benartzi2006a,toro2009ader,toro2015implicit}), has the same wave patterns as the Riemann solution of its corresponding homogeneous Euler equations. However, the singularity source term affects the wave patterns of the Riemann solution. In fact, as can be seen from \cite{2007Condensation}, the solution of Riemann problem (\ref{governing equations}) and (\ref{initial condition}) may be a four waves structure or a five waves structure, which is far different from the wave structure of the CRP.

In this work, an analytical frame, which is called the double CRPs frame, is proposed to construct the exact solution of the Riemann problem (\ref{governing equations}) and (\ref{initial condition}). It regards the fluids on both sides of the heating point as two separate CRPs, and gives the upper limit of the number of waves first. The solutions of the two CRPs are then coupled on the premise of maintaining the physical properties of the heating point. Depending on the heating properties and gasdynamics properties, the extra waves are deleted and the type of wave is finally determined. One advantage of the double CRPs frame is that it is independent of whether the fluid at the heating point is supersonic or subsonic. Under this frame we demonstrate three possible structures of the solution, which are verified by numerical tests in Section\ref{sec:algorithm}.

The text is arranged as follows. In Section\ref{sec:heating}, we will introduce the solution of the heating equations and derive several useful properties. In Section\ref{sec:structure}, an analytical frame of double CRPs coupling will be introduced to constructively solve the Riemann problem (\ref{governing equations}) and (\ref{initial condition}). Under this frame We will prove that there are at most three types of the solution. In Section\ref{sec:uniqueness}, the structure of the solution will be associated with the Mach number of the incoming flow to illustrate the uniqueness of the solution. In Section\ref{sec:algorithm}, we will give an iterative method for the solution of each structure. Then we will give five tests for the Riemann problem (\ref{governing equations}) and (\ref{initial condition}) with different initial conditions and heating coefficients. Finally, a brief summary will be given in Section\ref{sec:conclusions}.

\section{Solutions of the heating equations}
\label{sec:heating}
According to \cite{schnerr2005unsteadiness}, there are two branches of the solution to the heating equations (\ref{balance equations2}). We can choose the physical solution from these two branches by the following property.
\begin{Property}
	If $M_-<1$, then $M_+\leq 1$. If $M_->1$, then $M_+\geq 1$. 
\end{Property}

The physical solution corresponds to the weak solution in \cite{schnerr2005unsteadiness}. In this paper the expression of Dongen et al.\cite{2002On} is adopted. For the heat addition of subsonic flow, the solution is

\begin{subequations}
	\label{subsonic solution}
	\begin{align}
		&I\equiv \sqrt{\left( \gamma+\frac{1}{{M_-}^2}\right) ^2-2(\gamma+1)\left( \frac{1}{{M_-}^2}+\frac{\gamma-1}{2}\right) (1+\kappa)},\\
		&\frac{u_+}{u_-}=\frac{\rho_-}{\rho_+}=\frac{1}{\gamma+1}\left( \gamma+\frac{1}{{M_-}^2}- I\right), \\
		&\frac{p_+}{p_-}=\frac{{M_-}^2}{\gamma+1}\left( \gamma+\frac{1}{{M_-}^2}+ \gamma I\right), \\
		&M_+=\sqrt{\frac{\gamma+\frac{1}{{M_-}^2}- I}{\gamma+\frac{1}{{M_-}^2}+ \gamma I}}.
	\end{align}
\end{subequations}

For the heat addition of supersonic flow, the solution is

\begin{subequations}
	\label{supersonic solution}
	\begin{align}
		&I\equiv \sqrt{\left( \gamma+\frac{1}{{M_-}^2}\right) ^2-2(\gamma+1)\left( \frac{1}{{M_-}^2}+\frac{\gamma-1}{2}\right) (1+\kappa)},\\
		&\frac{u_+}{u_-}=\frac{\rho_-}{\rho_+}=\frac{1}{\gamma+1}\left( \gamma+\frac{1}{{M_-}^2}+ I\right), \\
		&\frac{p_+}{p_-}=\frac{{M_-}^2}{\gamma+1}\left( \gamma+\frac{1}{{M_-}^2}- \gamma I\right), \\
		&M_+=\sqrt{\frac{\gamma+\frac{1}{{M_-}^2}+ I}{\gamma+\frac{1}{{M_-}^2}- \gamma I}}.
	\end{align}
\end{subequations}

$M$ is Mach number. The advantage of this expression is that the ratios of variables before and after heat addition are only related to the upstream Mach number of the heating point. In order to make the solution reasonable, there is an upper bound on $k$ as follws.

\begin{equation}
	\label{maximum heat}
	k\leq k_{max}\mathop{=}\limits^{def}\frac{(1-{M_-}^2)^2}{2(\gamma+1){M_-}^2\left( 1+\frac{\gamma-1}{2}{M_-}^2\right)}.
\end{equation}

\begin{figure}[h]
	\label{figure maximum heat}
	\centering
	\includegraphics[height=15em]{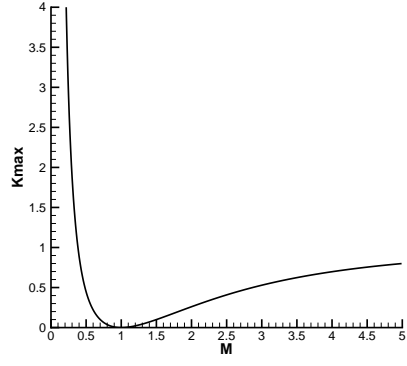}\par
	\caption{Relation between the maximum heating parameter and the upstream Mach number of the heating point.}
\end{figure}

$k_{max}$ is called maximum heating parameter. The relation between $k_{max}$ and ${M_-}$ is shown in Figure\ref{figure maximum heat}. When ${M_-}\to 0+$, $k_{max}$ approaches asymptotically $+\infty$. When ${M_-}\to \infty$, $k_{max}$ approaches asymptotically a finite number $\frac{1}{\gamma^2-1}$. The roots of the equation
\begin{equation}
	k=\frac{(1-{M}^2)^2}{2(\gamma+1){M}^2\left( 1+\frac{\gamma-1}{2}{M}^2\right)}
\end{equation}
are
\[
M^{1,2}=
\begin{cases}
	\sqrt{\frac{k(\gamma+1)+1\pm(\gamma+1)\sqrt{k(k+1)}}{1-k(\gamma^2-1)}}, \quad \text{if}\  k(\gamma^2-1)\ne 1\\
	\sqrt{\frac{1}{2[k(\gamma+1)+1]}}, \quad \text{if} \  k(\gamma^2-1)= 1	
\end{cases}
\]
We denote
\[
M_*=
\begin{cases}
	\sqrt{\frac{k(\gamma+1)+1-(\gamma+1)\sqrt{k(k+1)}}{1-k(\gamma^2-1)}},& \quad \text{if}\quad  k(\gamma^2-1)\ne 1\\
	\sqrt{\frac{1}{2[k(\gamma+1)+1]}},& \quad \text{if} \quad  k(\gamma^2-1)= 1	
\end{cases}
\]
\[
M_{**}=
\begin{cases}
	\sqrt{\frac{k(\gamma+1)+1+(\gamma+1)\sqrt{k(k+1)}}{1-k(\gamma^2-1)}},& \quad \text{if}\quad  k(\gamma^2-1)< 1\\
	+\infty,& \quad \text{if}\quad  k(\gamma^2-1)\geq 1	
\end{cases}
\]
It can easily be checked that $0<M_*<1$ and $M_{**}>1$, thus we arrive at the following conclusion.
\begin{Lemma}
	\label{lemma upwind condition}
	\mbox{}
	\begin{itemize}
		\item [(i)]
		If $k(\gamma^2-1)< 1$, then $M_-\leq M_*$ or $M_-\geq M_{**}$.
		\item [(ii)]
		If $k(\gamma^2-1)\geq 1$, then $M_-\leq M_*$.
		\item [(iii)]
		If $M_-=M_*$ or $M_-=M_{**}$, then $I=0$ and $M_+=1$.
	\end{itemize}
\end{Lemma}

\begin{Lemma}
	\label{lemma prandtl relation}
	If $k(\gamma^2-1)< 1$, we have
	\begin{equation}
		\label{prandtl relation}
		M_{*}=\sqrt{\frac{(\gamma-1)M_{**}^2+2}{2\gamma M_{**}^2-\gamma+1}},\quad
		M_{**}=\sqrt{\frac{(\gamma-1)M_*^2+2}{2\gamma M_*^2-\gamma+1}}.
	\end{equation}
\end{Lemma}
The proof is trivial.

The downstream fluid is called thermal choked if $M_+=1$. If the upstream fluid is sonic, then $k_{max}=0$, which means that any heat addition is not allowed for the sonic flow. We now turn to the relations of the upstream and downstream fluids.

\begin{Theorem}
	\label{theorem physical trend}
	For the heat additon of subsonic flow we have $u_+>u_-$, $p_+<p_-$, $\rho_+<\rho_-$ and $M_+>M_-$. For the heat additon of supersonic flow we have $u_+<u_-$, $p_+>p_-$, $\rho_+>\rho_-$ and $M_+<M_-$.
\end{Theorem}	

\begin{Proof}
	For the heat addition of subsonic fluid, $k>0$ implies
	\[ I< \sqrt{\left( \gamma+\frac{1}{{M_-}^2}\right) ^2-2(\gamma+1)\left( \frac{1}{{M_-}^2}+\frac{\gamma-1}{2}\right)}=\frac{1}{{M_-}^2}-1, \]
	\[ \frac{u_+}{u_-}> \frac{1}{\gamma+1}\left( \gamma+\frac{1}{{M_-}^2}+ 1-\frac{1}{{M_-}^2}\right)=1, \]
	\[ \frac{\rho_+}{\rho_-}=\frac{u_-}{u_+}<1, \]
	\[ \frac{p_+}{p_-}<\frac{{M_-}^2}{\gamma+1}\left( \gamma+\frac{1}{{M_-}^2}-\gamma+\frac{\gamma}{{M_-}^2}\right)=1, \]
	\[ M_+=\sqrt{\frac{\gamma+\frac{1}{{M_-}^2}- I}{\gamma+\frac{1}{{M_-}^2}+ \gamma I}}>\sqrt{\frac{\gamma+\frac{1}{{M_-}^2}- \frac{1}{{M_-}^2}+1}{\gamma+\frac{1}{{M_-}^2}+ \gamma (\frac{1}{{M_-}^2}-1)}}={M_-}.  \]
	For the heat addition of supsonic fluid, $k>0$ implies
	\[ I< \sqrt{\left( \gamma+\frac{1}{{M_-}^2}\right) ^2-2(\gamma+1)\left( \frac{1}{{M_-}^2}+\frac{\gamma-1}{2}\right)}=1-\frac{1}{{M_-}^2}, \]
	\[ \frac{u_+}{u_-}< \frac{1}{\gamma+1}\left( \gamma+\frac{1}{{M_-}^2}+ 1-\frac{1}{{M_-}^2}\right)=1, \]
	\[ \frac{\rho_+}{\rho_-}=\frac{u_-}{u_+}>1, \]
	\[ \frac{p_+}{p_-}>\frac{{M_-}^2}{\gamma+1}\left( \gamma+\frac{1}{{M_-}^2}- \gamma+\frac{\gamma}{{M_-}^2}\right)=1, \]
	\[ M_+=\sqrt{\frac{\gamma+\frac{1}{{M_-}^2}+ I}{\gamma+\frac{1}{{M_-}^2}- \gamma I}}<
	\sqrt{\frac{\gamma+\frac{1}{{M_-}^2}- \frac{1}{{M_-}^2}+1}{\gamma+\frac{1}{{M_-}^2}+ \gamma (\frac{1}{{M_-}^2}-1)}}={M_-}.  \]
\end{Proof}
Given any $k>0$ and $1<\gamma<3$, we define
\[
\begin{aligned}
	\phi(M)=\frac{1}{\gamma+1}\left( \gamma+\frac{1}{{M}^2}- I\right) ,\\
	\psi(M)=\frac{{M}^2}{\gamma+1}\left( \gamma+\frac{1}{{M}^2}+ \gamma I\right) ,
\end{aligned}
\]
where
\[
I=\sqrt{\left( \gamma+\frac{1}{{M}^2}\right) ^2-2(\gamma+1)\left( \frac{1}{{M}^2}+\frac{\gamma-1}{2}\right) (1+\kappa)}.
\]
According to \ref{subsonic solution} and \ref{supersonic solution}, $\frac{u_+}{u_-}$ and $\frac{p_+}{p_-}$ are both functions of $M_-$. For the subsonic heat addition, we have
\begin{equation}
	\label{subsonic solution function}
	\frac{u_+}{u_-}=\frac{\rho_-}{\rho_+}=\phi(M_-),\quad 
	\frac{p_+}{p_-}=\psi(M_-).
\end{equation}

A tedious compution gives the following theorem.
\begin{Theorem}
	\label{theorem subsonic derivation}
	$\phi'(M)>0,\quad \psi'(M)<0$.	
\end{Theorem}

\section{Structure of solution}
\label{sec:structure}
A basic assumption of our double CRPs frame is that $U_-(t)$ and $U_+(t)$ are both constant vectors. The solution $U(x,t)$ of the Riemann problem (\ref{governing equations}) and (\ref{initial condition}) satisfies the homogeneous Euler equations in both the left half $\{(x,t)|x<0, t\ge 0\}$ and the right $\{(x,t)|x>0, t\ge 0\}$. Consequently, $U(x,t)$ is self-similar. The exact solution in this paper refers to the self-similar solution under this frame.

\begin{figure}[h]
	\label{figure general structure}
	\centering
	\includegraphics[height=12em]{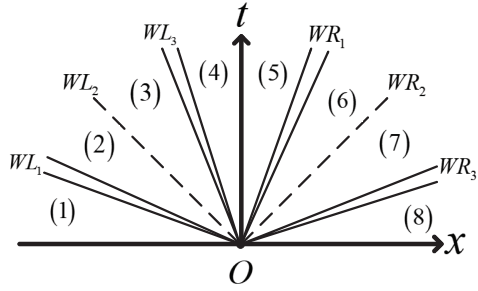}\par
	\caption{All possible waves for the Riemann problem (\ref{governing equations}) and (\ref{initial condition})}
\end{figure}

The exact solution $U(x,t)$ in the left half $\{(x,t)|x<0, t\ge 0\}$ satisfies the classical Euler equations, hence is the left half of the CRP solution with $U_1$ and $U(0-,t)$ as the left and right initial conditions. From the theory of CRP solution, we know that $U(x,t)$ in the left half consists three discontinuities at most, which are two genuinely nonlinear waves, namely shock waves or rarefaction waves, and a contact discontinuity corresponding to the characteristic fields $u-a$, $u+a$ and $u$, respectively. Similarly, $U(x,t)$ in the right half consists two genuinely nonlinear waves and a contact discontinuity at most.
\begin{Theorem}
	The self-similar solution of the Riemann problem (\ref{governing equations}) and (\ref{initial condition}) consists of seven discontinuities at most. They are a heating discontinuity at $x=0$, two genuinely nonlinear waves and a contact discontinuity left to $x=0$, two genuinely nonlinear waves and a contact discontinuity right to $x=0$ respectively, as shown in Figure\ref{figure general structure}.
\end{Theorem}

The elementary waves on the left and right sides are denoted as $WL_1$, $WL_2$, $WL_3$ and $WR_1$, $WR_2$, $WR_3$, respectively. The eight constant regions are labeled $(1)\sim (8)$, respectively. $U_1\sim U_8$ are the states in regoins $(1)\sim (8)$ and it is clear that $U_1=U_8$. Note that the t-axis is a discontinuity and $U_4\ne U_5$. We adopt similiar way to express the solution structure as in \cite{thanh2009the}. For examples, $S(U_1,U_2)$ and $R(U_1,U_2)$ mean two states $U_1$ and $U_2$ are connected by a shock and a rarefaction wave respectively, $C(U_2,U_3)$ means $U_2$ and $U_3$ are connected by a contact discontinuity, and $H(U_4,U_5)$ means $U_4$ and $U_5$ are connected by a heating discontinuity. All the six elementary waves in Figure\ref{figure general structure} can not exist at the same time. The next step in our double CRPs coupling method is to eliminate the redundant waves according to the heat addition properties and the gasdynamics properties. We will give the main results in Theorem\ref{theorem solution structure}.

\begin{Theorem}
	\label{theorem solution structure}
	Under the double CRPs frame, there are three different structures of the exact solution to the Riemann problem (\ref{governing equations}) and (\ref{initial condition}) as follows.
	\begin{itemize}
		\item [(1)]
		Type 1: if $M_4<M_5<1$, then the structure could be \[S(U_1,U_4)\oplus H(U_4,U_5)\oplus C(U_5,U_7)\oplus S(U_7,U_8)\].
		\item [(2)]
		Type 2: if $M_4<M_5=1$, then the structure could be \[S(U_1,U_4)\oplus H(U_4,U_5)\oplus R(U_5,U_6)\oplus C(U_6,U_7)\oplus S(U_7,U_8)\].
		\item [(3)]
		Type 3: if $M_4>M_5\geq 1$, then the structure could be \[H(U_1,U_5)\oplus R(U_5,U_6)\oplus C(U_6,U_7)\oplus S(U_7,U_8)\].
	\end{itemize}
\end{Theorem}

\begin{figure}[h]
	\label{figure preliminary structure}
	\centering
	\includegraphics[height=10em]{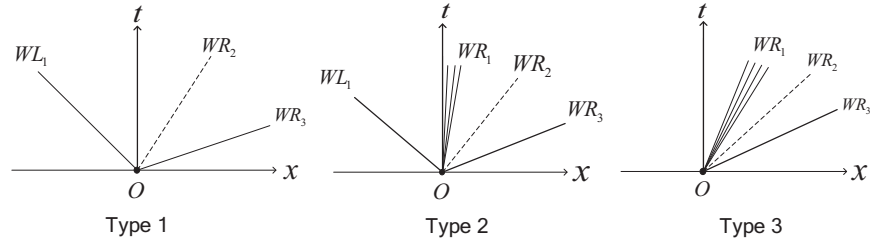}\par
	\caption{Three different structures of the exact solution to the Riemann problem (\ref{governing equations}) and (\ref{initial condition}).}
\end{figure}

Figure\ref{figure preliminary structure} depicts the three wave patterns in Theorem\ref{theorem solution structure}. They will be denoted by Type 1, Type 2 and Type 3, respectively. According to the Mach number at the heating point and whether the thermal choked state appears, the proof of the Theorem\ref{theorem solution structure} falls naturally into three parts: Lemma\ref{lemma structure 1}, Lemma\ref{lemma structure 2} and Lemma\ref{lemma structure 3}, which correspond to Type 1, Type 2 and Type 3, respectively.

\begin{Lemma}
	\label{lemma structure 1}
	Under the double CRPs frame, if $M_4<M_5<1$, then $WL_2$, $WL_3$ and $WR_1$ do not exist, and then $u_5>u_4>0$; if $WL_1$ and $WR_3$ exist, they are both shock waves.
\end{Lemma}

\begin{Proof}
	We first outline the proof. The problem is reduced to the determination of the upstream fluid $U_4$ and the downstream fluid $U_5$ (either of which can be determined by the other) at the heating point, such that $U_4$ and $U_1$ can be connected by the fundamental waves with non-positive speed, and that $U_5$ and $U_8=U_1$ can be connected by the fundamental waves with non-negative speed.
	
	According to (\ref{balance equation}), it follows that $u_4$ and $u_5$ have the same sign. $WL_3$ and $WR_1$ do not exist since $M_4<1$ and $M_5<1$. There are three cases for the sign of the velocity at  heating point, as follows.
	\begin{itemize}
		\item [(1)]
		If $u_4>0$, then $u_5>0$ and $WL_2$ does not exist;
		\item [(2)]
		If $u_4<0$, then $u_5<0$ and $WR_2$ does not exist;
		\item [(3)]
		If $u_4=u_5=0$, then both the speeds of $WL_2$ and $WR_2$ are zero and the widths of region $3$ and region $6$ are both zero.
	\end{itemize}
	
	These three structures are shown in Figure\ref{figure temporary sructure}.
	
	\begin{figure}[h]
		\label{figure temporary sructure}
		\centering
		\includegraphics[height=10em]{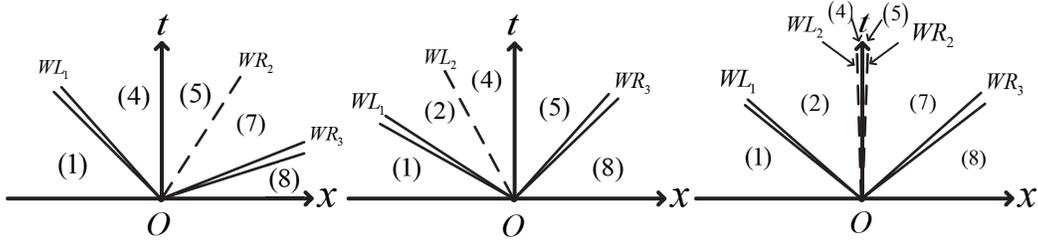}\par
		\caption{Three possible structures when subsonic flow is heated without thermal choke.}
	\end{figure}
	
	If $u_4<0,u_5<0$, then $WL_1$ is a shock wave and $WR_3$ is a rarefaction wave. From Theorem\ref{theorem physical trend} we have $p_4<p_5$. The relations of the elenmentary wave (see for instance \cite{0Riemann}) imply $p_1\leq p_2=p_4$ and $p_5\leq p_8$. It results $p_5\leq p_1$, which is contradictory with $p_5>p_4\geq p_1$, hence the structure of $u_4<0,u_5<0$ does not exist.
	
	Using a similar method we can prove that the structure of $u_4=u_5=0$ dose not exist. The structure of the solution is the left figure in Figure 5.
	
	The next thing to do is to prove that $WL_1$ is a shock. If the assertion would not hold, then $WL_1$ is a rarefaction wave, which implies $p_7=p_5<p_4\leq p_1=p_8$. It follows that $WR_3$ is a rarefaction wave and $u_7=u_5>u_4\geq u_1$. Applying the conservations of mass and momentum at the control volume $[x_0,x_4]\times [0,T]$ in the $x-t$ space, as shown in the left figure of Figure\ref{figure integral region}, we have
	\[
	\begin{aligned}
		\int_{x_0}^{x_4}\rho(x,T) dx&=\int_{x_0}^{x_4}\rho(x,0) dx,\\
		\int_{x_0}^{x_4}\rho(x,T)u(x,T) dx&=\int_{x_0}^{x_4}\rho(x,0)u(x,0) dx.
	\end{aligned}
	\]
	$U$ is constant at each region, thus
	\begin{equation}
		\label{lemma proof formula1}
		\begin{aligned}
			&\int_{x_0}^{x_1}\rho(x,T)dx+\rho_4|x_1|+\rho_5|x_2|+\rho_7|x_3-x_2|+\int_{x_3}^{x_4}\rho(x,T)dx\\
			&=\int_{x_0}^{x_1}\rho(x,T)\frac{u(x,T)}{u_1}dx+\rho_4 |x_1|\frac{u_4}{u_1}+\rho_5|x_2|\frac{u_5}{u_1}+\rho_7|x_3-x_2|\frac{u_7}{u_1}
			+\int_{x_3}^{x_4}\rho(x,T)\frac{u(x,T)}{u_1}dx.
		\end{aligned}
	\end{equation}
	The velocity inside the rarefaction wave is monotonous, it follows that $u(x,T)\geq max\{u_1,u_4\}\geq u_1$ for $x_0\leq x\leq x_1$ and $u(x,T)\geq max\{u_5,u_8\}=u_5>u_1$ for $x_3\leq x\leq x_4$, which imply
	\[\int_{x_0}^{x_1}\rho(x,T)\frac{u(x,T)}{u_1}dx\geq \int_{x_0}^{x_1}\rho(x,T)dx,\] 
	\[\int_{x_3}^{x_4}\rho(x,T)\frac{u(x,T)}{u_1}dx\geq \int_{x_3}^{x_4}\rho(x,T)dx.\] 
	$\frac{u_4}{u_1}\geq 1$, $\frac{u_5}{u_1}> 1$ and $\frac{u_7}{u_1}> 1$ hold for the right side of Lemma\ref{lemma proof formula1}. Substituting the above inequalities to Lemma\ref{lemma proof formula1}, we have $x_2=x_3=0$ and $u_5=0$, which leads to a contradiction.
	
	Finally, we have to show that $WR_3$ is a shock. If the assertion is false, then $WR_3$ is a non-degenerate rarefaction wave, which implies $u_4< u_5=u_7\leq u_8=u_1$. Applying the conservations of mass and momentum at the control volume $[x_0,x_3]\times [0,T]$ in the $x-t$ space, as shown in the right figure of Figure\ref{figure integral region}, we have
	
	\begin{figure}[h]
		\label{figure integral region}
		\centering
		\includegraphics[height=10em]{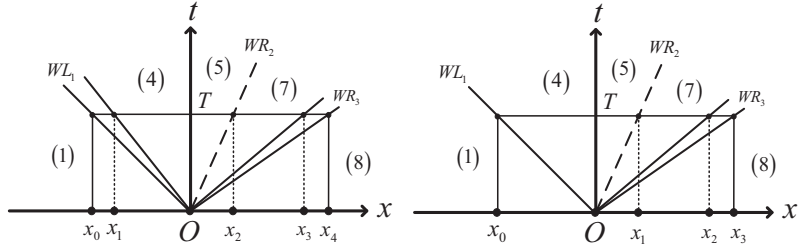}\par
		\caption{The integral diagram in the proof of Lemma\ref{lemma structure 1}. Left: $WL_1$ is a rarefaction wave and $WR_3$ is a rarefaction wave. Right: $WL_1$ is a shock and $WR_3$ is a rarefaction wave.}
	\end{figure}
	
	\[
	\begin{aligned}
		\int_{x_0}^{x_3}\rho(x,T) dx&=\int_{x_0}^{x_3}\rho(x,0) dx,\\
		\int_{x_0}^{x_3}\rho(x,T)u(x,T) dx&=\int_{x_0}^{x_3}\rho(x,0)u(x,0) dx.
	\end{aligned}
	\]
	Thus
	\[
	\begin{aligned}
		&\rho_4|x_0|+\rho_5|x_1|+\rho_7|x_2-x_1|+\int_{x_2}^{x_3}\rho(x,T)dx\\
		&=\rho_4 |x_0|\frac{u_4}{u_1}+\rho_5|x_1|\frac{u_5}{u_1}+\rho_7|x_2-x_1|\frac{u_7}{u_1}
		+\int_{x_2}^{x_3}\rho(x,T)\frac{u(x,T)}{u_1}dx.
	\end{aligned}
	\]
	For the velocity $u$ insides the rarefaction wave $WR_3$, we have $u(x,T)\leq max(u_5,u_8)=u_8=u_1$, which implies
	\[\int_{x_2}^{x_3}\rho(x,T)\frac{u(x,T)}{u_1}dx\leq \int_{x_2}^{x_3}\rho(x,T)dx.\] 
	$\frac{u_4}{u_1}<1$, $\frac{u_5}{u_1}\leq 1$ and $\frac{u_7}{u_1}\leq 1$ hold for the right side of Lemma\ref{lemma proof formula1}. We see at once that $\frac{u_7}{u_1}=1$ and $U_7=U_8$, which lesds to a contradiction.
\end{Proof}

The following two lemmas can be proved using a similar method.

\begin{Lemma}
	\label{lemma structure 2}
	Under the double CRPs frame, if $M_4<M_5=1$, then $WL_2$ and $WL_3$ do not exist, and then $u_5>u_4>0$; if $WL_1$, $WR_1$ and $WR_3$ exist, they are a shock, a rarefaction wave and a shock, respectively.
\end{Lemma}

\begin{Lemma}
	\label{lemma structure 3}
	Under the double CRPs frame, if $M_4>M_5\ge 1$, then $WL_1, $$WL_2$ and $WL_3$ do not exist, and then $u_5>u_4>0$; if $WR_1$ and $WR_3$ exist, they are a rarefaction wave and a shock, respectively.
\end{Lemma}

Compared with Type 1 and Type 3, Type 2 has one more wave, and the condition $M_5=1$ is used to match the number of conditions at this time. In the remaindeer of this paper, each wave and each constant region are denoted as in Figure\ref{figure preliminary structure}.

\section{Uniqueness of solution}
\label{sec:uniqueness}

In this section, we will make a preliminary exploration of the uniqueness of the solution.

A more refined analysis of the exact solution needs to involve the quantitative relations of the fundamental waves and the heating discontinuity. We begin by defining five functions to connect the states on the left and right sides of the shock wave $WL_1$, which appears in the solution of Type 1 and Type 2.
\[
\begin{aligned}
	&f_1(M_1,M_2)\mathop{=}\limits^{def} \frac{(\gamma+1)(M_1-M_2)^2}{(\gamma-1)(M_1-M_2)^2+2},\\
	&f_2(M_1,M_2)\mathop{=}\limits^{def}\frac{2\gamma (M_1-M_2)^2-\gamma+1}{\gamma+1},\\
	&f_3(M_1,M_2)\mathop{=}\limits^{def}\frac{(\gamma-1)(M_1-M_2)^2+2}{(\gamma+1)(M_1-M_2)^2},\\
	&f_4(M_1,M_2)\mathop{=}\limits^{def}f_3(M_1,M_2)+\frac{M_2}{M_1}(1-f_3(M_1,M_2)),\\
	&f_5(M_1,M_2)\mathop{=}\limits^{def}\frac{((\gamma-1)M_1+2M_2)(M_1-M_2)+2}{[2\gamma(\gamma-1)(M_1-M_2)^4+(6\gamma-\gamma^2-1)(M_1-M_2)^2-2(\gamma-1)]^{\frac{1}{2}}}.
\end{aligned}
\]
Applying the Rankine-Hugoniot conditions
\[
F(U_4)-F(U_1)=s_L(U_4-U_1),
\]
where $s_L$ is the speed of $WL_1$, we deduce that
\begin{equation}
	\label{Lshock formula}
	\begin{aligned}
		&{\rho_4}/{\rho_1}=f_1(M_1,M_{SL}),\quad
		&{p_4}/{p_1}=f_2(M_1,M_{SL}),\\
		&({u_4-s_L})/({u_1-s_L})=f_3(M_1,M_{SL}),\quad 
		&{u_4}/{u_1}=f_4(M_1,M_{SL}),\\
		&{a_4}/{a_1}=f_2(M_1,M_{SL})^{\frac{1}{2}}f_1(M_1,M_{SL})^{-\frac{1}{2}},\quad
		&M_4=f_5(M_1,M_{SL}),
	\end{aligned}
\end{equation}
where $M_{SL}=s_L/a_1$ is called the shock Mach number of $WL_1$.

We now turn to the overall relations of the solution. According to (\ref{subsonic solution function}), the relations of the stationary discontinuity at the origin are
\begin{subequations}
	\label{heat addition formula}
	\begin{align}
		\frac{p_5}{p_4}=\psi(M_4)=\psi(f_5(M_1,M_{SL})),\\
		\frac{u_5}{u_4}=\phi(M_4)=\phi(f_5(M_1,M_{SL})).
	\end{align}
\end{subequations}

It is obtained from the shock relation of $WR_3$ that
\[ u_5-u_8=\sqrt{\frac{\beta p_8}{\rho_8}}\frac{p_5/p_8-1}{\sqrt{1+\tau p_5/p_8}} \]
where $\beta=2/(\gamma-1),\tau=(\gamma+1)/(\gamma-1)$.

Substituting (\ref{Lshock formula}) and (\ref{heat addition formula}) into the above equation, we get
\begin{equation}
	\label{left equation}
	M_1(f_4(M_1,M_{SL})\phi(M_4)-1)
	-\sqrt{\frac{\beta}{\gamma}}\frac{f_2(M_1,M_{SL})\psi(M_4)-1}{\sqrt{1+\tau f_2(M_1,M_{SL})\psi(M_4)}}=0.
\end{equation}

The equation (\ref{Lshock formula}) forms a system of equations for $M_{SL}$ and $M_4$ as follows.
\begin{equation}
	\label{left equations}
	\left\lbrace 
	\begin{aligned}
		&M_1(f_4(M_1,M_{SL})\phi(M_4)-1)
		-\sqrt{\frac{\beta}{\gamma}}\frac{f_2(M_1,M_{SL})\psi(M_4)-1}{\sqrt{1+\tau f_2(M_1,M_{SL})\psi(M_4)}}=0\\
		&M4=\frac{((\gamma-1)M_1+2M_{SL})(M_1-M_{SL})+2}{\sqrt{2\gamma(\gamma-1)(M_1-M_{SL})^4+(6\gamma-\gamma^2-1)(M_1-M_{SL})^2-2(\gamma-1)}}
	\end{aligned}
	\right. 
\end{equation}

Another preliminary theorem is to conclude that $M_{SL}$ can be expressed as a function of $M_1$ and $M_4$ in the solution of Type 1 or Type 2. From the Lax entropy condition of the shock (see \cite{Peter1973Hyperbolic}), we have
\begin{equation}
	\label{lax condition}
	M_1-M_{SL}=\frac{u_1-s_L}{a_1}\ge 1.
\end{equation}

\begin{Lemma}
	\label{lemma solvability}
	Given any $M_1$ and $M_4$, the $M_{SL}$ satisfying $M_4=f_5(M_1, M_{SL})$ and $M_1-M_{SL}\geq 1$ exists and is unique.
\end{Lemma}

\begin{Proof}
	The formula $M_4=f_5(M_1, M_{SL})$ has another form.
	\[
	\begin{aligned}
		&M_4\sqrt{(\gamma-1)(M_1-M_{SL})^2+2}\sqrt{2\gamma (M_1-M_{SL})^2-\gamma+1}\\
		=&-2(M_1-M_{SL})^2+(\gamma+1)M_1(M_1-M_{SL})+2.
	\end{aligned}
	\]
	Dividing both sides of the above equation by $M_1-M_{SL}$, we have
	\[
	\begin{aligned}
		&M_4\sqrt{(\gamma-1)(M_1-M_{SL})^2+2}\sqrt{2\gamma-\frac{\gamma+1}{ (M_1-M_{SL})^2}}\\
		=&-2(M_1-M_{SL})+(\gamma+1)M_1(M_1-M_{SL})+\frac{2}{(M_1-M_{SL})^2}.
	\end{aligned}
	\]	
	We define
	\[
	\begin{aligned}
		&\sigma_1(x)\mathop{=}\limits^{def}M_4\sqrt{(\gamma-1)x^2+2}\sqrt{2\gamma-\frac{\gamma+1}{x^2}},\\
		&\sigma_2(x)\mathop{=}\limits^{def}2x-(\gamma+1)M_1-\frac{2}{x},\\
		&\sigma(x)\mathop{=}\limits^{def}\sigma_1(x)+\sigma_2(x).
	\end{aligned}
	\]
	The domains of $\sigma_1(x)$, $\sigma_2(x)$ and $\sigma(x)$ are $\{x|x\ge 1\}$. Both $\sigma_1(x)$ and $\sigma_2(x)$ are monotone increasing functions, hence $\sigma(x)$ is a monotone increasing function. It is clear that
	\[
	\sigma(1)=(\gamma+1)(M_4-M_1)\le 0,
	\]
	\[
	\lim_{x\to \infty}\sigma(x)=\lim_{x\to \infty}\sigma_1(x)+\lim_{x\to \infty}\sigma_2(x)=\infty.
	\]
	Consequently, the root of $\sigma(x)$ exists and is unique. $M_1-M_{SL}$ is the root of $\sigma(x)$, and this completes the proof.
\end{Proof}

For the solutions of Type 1 and Type 2, it clear that
\[ M_4=\frac{((\gamma-1)M_1+2M_{SL})(M_1-M_{SL})+2}{\sqrt{2\gamma(\gamma-1)(M_1-M_{SL})^4+(6\gamma-\gamma^2-1)(M_1-M_{SL})^2-2(\gamma-1)}}.\]
Treating $M_1-M_{SL}$ as the variable, we have
\begin{equation}
	\label{polynomial function}
	A(M_1-M_{SL})^4+B(M_1-M_{SL})^3+C(M_1-M_{SL})^2+D(M_1-M_{SL})+E=0,
\end{equation}
where \[
\begin{aligned}
	&A=4-2\gamma(\gamma-1)M_4^2,\\
	&B=-4(\gamma+1)M_1,\\
	&C=(\gamma+1)^2M_1^2-8-(6\gamma-\gamma^2-1)M_4^2,\\
	&D=4(\gamma+1)M_1,\\
	&E=4+2(\gamma-1)M_4^2.
\end{aligned}
\]
$M_1-M_{SL}$ is the root of a fourth-order polynomial function (\ref{polynomial function}). According to Lemma\ref{lemma solvability}, We have proved the following theorem.
\begin{Theorem}
	\label{theorem expression}
	For the solutions of Type 1 and Type 2, $M_{SL}$ can be expressed as $M_{SL}=f_6(M_1,M_4)$, which satisfies $M_4=f_5(M_1,M_{SL})$ and $M_1-M_{SL}\geq 1$.
\end{Theorem}

Now we turn to the uniqueness of the solution. The proof consists of two parts, one involes the states at the intermediate regions (Lemma\ref{lemma uniqueness1}), and the other is a division of these three structures (Lemma\ref{lemma uniqueness2}). Both of these are related to the Mach number $M_1$ of the initial flow.

We describe the strength of each nonlinear wave in terms of the ratio of pressures on both sides of this wave.
\begin{Lemma}
	\label{lemma uniqueness1}
	For these three types of solutions, the strength of each nonlinear wave is determined by the Mach number $M_1$, independent of other variables.
\end{Lemma}
A proof of this lemma will be given in Appendix\ref{proof uniqueness1}.

We define
\begin{equation}
	\label{function X}
	X(M_1,M_4)\mathop{=}\limits^{def}M_1(f_4(M_1,M_{SL})\phi(M_4)-1)
	-\sqrt{\frac{\beta}{\gamma}}\frac{f_2(M_1,M_{SL})\psi(M_4)-1}{\sqrt{1+\tau f_2(M_1,M_{SL})\psi(M_4)}},
\end{equation}
and
\begin{equation}
	\label{function Y}
	Y(M_1)\mathop{=}\limits^{def}X(M_1,M_*),
\end{equation}
where $ M_{SL}=f_6(M_1,M_4)$.

If $k(\gamma^2-1)\geq 1$, then $M_{**}$ does not exist. Therefore the upstream flow of the heating point must be a subsonic flow. At this time, there are only two possible structures: Type 1 and Type 2. In this case, the downstream fluid must be subsonic. When $k(\gamma^2-1)< 1$, all three structures are possible. We associate the three structures with $M_1$ by the following lemma. 

\begin{Lemma}
	\label{lemma uniqueness2}
	Under the double CRPs frame, the structures of the exact solution can only be Type 1 or Type 2 for $k(\gamma^2-1)\geq 1$, and all three structures are possible for $k(\gamma^2-1)< 1$. This type of solution satisfies the following conditions.
	\begin{itemize}
		\item [(i)]
		For the solution of Type 1, $Y(M_1)\geq 0$ holds;
		\item [(ii)]
		For the solution of Type 2, $Y(M_1)\leq 0$ and $M_1 \geq M_*$ hold, and $M_1\leq M_{**}$ holds at the case of $k(\gamma^2-1)< 1$;
		\item [(iii)]
		For the solution of Type 3, $M_1\geq M_{**}$ holds.
	\end{itemize}
\end{Lemma}
A proof of this lemma will be given in Appendix\ref{proof uniqueness2}.

It can be seen from Lemma\ref{lemma uniqueness2} that the structure of $M_1$ equaling to the root of $Y(M_1)$ is the demarcation structure of Type 1 and Type 2, which is
\[S(U_1,U_4)\oplus H(U_4,U_5)\oplus C(U_5,U_7)\oplus S(U_7,U_8)\quad \text{with}\quad M_4=M_*.\] 
And the structure of $M_1$ equaling to $M_{**}$ is the demarcation structure of Type 2 and Type 3, which is
\begin{equation}
	\label{limit structure 1}
	S(U_1,U_4)\oplus H(U_4,U_5)\oplus R(U_5,U_6)\oplus C(U_6,U_7)\oplus S(U_7,U_8) \quad \text{with}\quad s_L=0,
\end{equation}
or
\begin{equation}
	\label{limit structure 2}
	H(U_1,U_5)\oplus R(U_5,U_6)\oplus C(U_6,U_7)\oplus S(U_7,U_8) \quad \text{with}\quad M_5=1.
\end{equation}

\begin{Remark}
	The structure of (\ref{limit structure 1}) is a limit structure of Type 2. $WL_1$ is a normal shock and $M_4=M_*$ hold for this structure, and it is clear that $M_1=M_{**}$. For the structure of (\ref{limit structure 2}), $M_5=1$ implies $M_1=M_{**}$. Therefore the structure of (\ref{limit structure 1}) equals to the structure of (\ref{limit structure 2}).
\end{Remark}

According to Lemma\ref{lemma uniqueness1} and Lemma\ref{lemma uniqueness2}, we can establish the following theorem.
\begin{Theorem}	
	\label{theorem uniqueness}
	Under the double CRPs frame, we have
	\begin{itemize}
		\item [(i)]
		if $k(\gamma^2-1)\geq 1$, the solution of Riemann problem (\ref{governing equations}) and (\ref{initial condition}) is unique;
		\item [(ii)]
		if $k(\gamma^2-1)< 1$, the solution of Riemann problem (\ref{governing equations}) and (\ref{initial condition}) is unique under the assumption that the root of $Y(M_1)$ is not greater than $M_{**}$.
	\end{itemize}
\end{Theorem}

If the root of $Y(M_1)$ is greater than $M_{**}$, the uniqueness of the solution has not been proved. To compare the size of the root of $Y(M1)$ and $M_{**}$, we define
\[
R(\gamma,k)=\{M|Y(M)=0\},\quad 
T(\gamma,k)=M_{**}-R(\gamma,k).
\]

\begin{figure}[h]
	\label{figure uniqueness test}
	\centering
	\includegraphics[height=6.5cm]{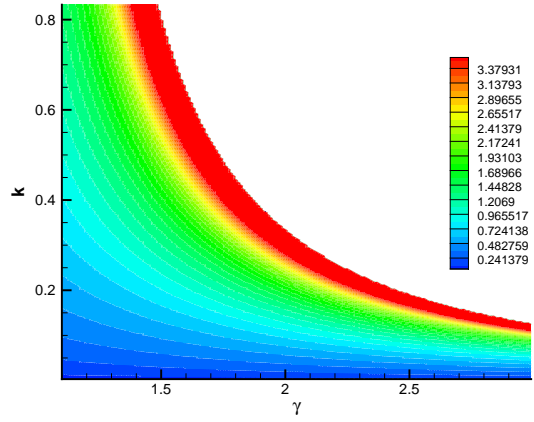}\par
	\caption{The contour of function $T(\gamma,k)$}
\end{figure}

Figure\ref{figure uniqueness test} shows the contour of the function $T$. The curved boundary on the upper right is the curve of $k(\gamma^2-1)=1$. From Figure\ref{figure uniqueness test} it can be found $T(\gamma,k)>0$, therefore the root of $Y(M_1)$ is smaller than $M_{**}$ and the assumption in the second part of Theorem\ref{theorem uniqueness} is true. Thus we give the following conjecture.
\begin{Conjecture}
	\label{conjecture}
	Under the double CRPs frame, the solution of Riemann problem (\ref{governing equations}) and (\ref{initial condition}) is unique for any given $\gamma$ and $\kappa$, and the structure of the self-similar solution is determined by $M_1$, as shown below.
	\begin{itemize}
		\item [(i)]
		If $Y(M_1)\geq 0$, the structure of the solution is Type 1.
		\item [(ii)]
		If $Y(M_1)\geq 0$ and $M_1\leq M_{**}(k(\gamma^2-1)<1)$, the structure of the solution is Type 2.
		\item [(iii)]
		If $k(\gamma^2-1)<1$ and $M_1\geq M_{**}$, the structure of the solution is Type 3.
	\end{itemize}
\end{Conjecture}

\section{Algorithm of Solution and Verification}
\label{sec:algorithm}

\begin{algorithm}[htbp]
	 \small
	\caption{The Construction Algorithm of the Exact Solution}  
	\label{algorithm}  
	\begin{algorithmic} 
		\Require  
		the initial condition:  $\rho_1$,   $u_1$,   $p_1$,  $M_1$, $a_1$; 
		\Ensure  
		the intermediate states: $W_4=(\rho_4, u_4, p_4)$,  $W_5=(\rho_5, u_5, p_5)$,  $W_6=(\rho_6, u_6, p_6)$,  $W_7=(\rho_7, u_7, p_7)$;
		\State compute $M_*$ and $M_{**}$;
		\State compute $Y(M_1)$;
		\If {$M_1\geq M_{**}$}  
		\State   $W_4\leftarrow W_1$;
		\State   compute $W_5$ by (\ref{supersonic solution});
		\ElsIf{$Y(M_1)\geq 0$}
		\State set the initial values of the bisection method: $x_1$, $x_2$;
		\State $x_0\leftarrow (x_1+x_2)/2$;
		\State set the value of termination condition: $error$;
		\While {$|x_1-x_2|>error$}
		\State $f_1\leftarrow ((\gamma-1)M_1+2x_1)(M_1-x_1)+2-M_*(2\gamma(\gamma-1)(M_1-x_1)^4+(6\gamma -\gamma^2-1)(M_1-x_1)^2-2(\gamma-1))^{\frac{1}{2}}$;
		\State $f_2\leftarrow ((\gamma-1)M_1+2x_2)(M_1-x_1)+2-M_*(2\gamma(\gamma-1)(M_1-x_2)^4+(6\gamma -\gamma^2-1)(M_1-x_2)^2-2(\gamma-1))^{\frac{1}{2}}$;
		\If {$f_1\times f_2<0$}
		\State   $x_0\leftarrow x_1$;
		\Else
		\State   $x_1\leftarrow x_0$;
		\EndIf
		\State $x_0\leftarrow (x_1+x_2)/2$;
		\EndWhile
		\State   $\rho_4\leftarrow \rho_1f_1(M_1,x_0)$, $u_4\leftarrow u_1f_4(M_1,x_0)$, $p_4\leftarrow p_1f_2(M_1,x_0)$, $M_4\leftarrow M_*$;
		\State   compute $W_5$ by (\ref{subsonic solution});
		\State   $p_5\leftarrow p_4\psi(M_*)$, $u_5\leftarrow u_4\phi(M_*)$, $\rho_5\leftarrow (\rho_4u_4)/u_5$;
		\Else  \Statex\ //$f_L$ and $f_R$ are given by (\ref{shock pressure formula}).	
		\State set the initial values of the bisection method: $y_1$, $y_2$;
		\State $y_0\leftarrow (y_1+y_2)/2$;
		\State set the value of termination condition: $error$;
		\While {$|y_1-y_2|>error$}
		\State $m_1\leftarrow (u_1-f_L(y_1))\left[\frac{\rho_1}{\gamma}\frac{\gamma+1+(\gamma-1)p_1/y_1}{(\gamma-1)y_1+(\gamma+1)p_1} \right] ^{\frac{1}{2}}$, 
		$m_2\leftarrow (u_1-f_L(y_2))\left[\frac{\rho_1}{\gamma}\frac{\gamma+1+(\gamma-1)p_1/y_2}{(\gamma-1)y_2+(\gamma+1)p_1} \right] ^{\frac{1}{2}}$;
		\State $g_1\leftarrow u_1+f_R(y_1\psi(m_1))-\phi(m_1)(u_1-f_L(y_1))$, 
		$g_2\leftarrow u_1+f_R(y_2\psi(m_2))-\phi(m_1)(u_1-f_L(y_2))$;
		\If {$g_1\times g_2<0$}
		\State   $y_0\leftarrow y_1$;
		\Else
		\State   $y_1\leftarrow y_0$;
		\EndIf
		\State $y_0\leftarrow (y_1+y_2)/2$;
		\EndWhile
		\State   $\rho_4\leftarrow \rho_1\frac{(\gamma-1)p_1+(\gamma+1)y_0}{(\gamma-1)y_0+(\gamma+1)p_1}$, $u_4\leftarrow u_1-f_L(y_0)$, $p_4\leftarrow y_0$;
		\State $M_1\leftarrow (u_1-f_L(y_0))\left[\frac{\rho_1}{\gamma}\frac{\gamma+1+(\gamma-1)p_1/y_0}{(\gamma-1)y_0+(\gamma+1)p_1} \right] ^{\frac{1}{2}}$;
		\State   compute $W_5$ by (\ref{subsonic solution});
		\EndIf 
		\State compute $W_6$ and $W_7$ by a Riemann solver $CRP(W_5,W_8)$. 
	\end{algorithmic}  
\end{algorithm}

Lemma\ref{lemma uniqueness2} makes it legitimate to apply $M_1$ to determine the structure of the exact solution. Once the structure is determined, the solution becomes easy. A algorithm for the solution is given in Algorithm\ref{algorithm}.

In the algorithm of Type 2, we use the bisection method to solve 
\begin{equation}
	\label{structure2 mach number}
	M_*=\frac{((\gamma-1)M_1+2M_{SL})(M_1-M_{SL})+2}{\sqrt{2\gamma(\gamma-1)(M_1-M_{SL})^4+(6\gamma-\gamma^2-1)(M_1-M_{SL})^2-2(\gamma-1)}}.
\end{equation}
In the algorithm of Type 1, we do not iteratively solve \ref{left equations}, but use the pressure equations
\[
u_1+f_R(p_4\psi(M_4(p_4)))-\phi(M_4(p_4))(u_1-f_L(p_4))=0,
\]
where
\[
M_4(p_4)=[u_1-f_L(p_1)]\left[ \frac{\rho_1}{\gamma}\frac{\gamma+1+(\gamma-1)p_1/p_4}{(\gamma-1)p_4+(\gamma+1)p_1}\right] ^{\frac{1}{2}},
\]
\begin{equation}
	\label{shock pressure formula}
	\begin{aligned}
		f_L(p)=
		\begin{cases}
			(p-p_1)\left[ \frac{A_L}{p+B_L}\right] ^{\frac{1}{2}}, \quad if \ p>p_1\\
			\frac{2a_1}{\gamma-1}\left[ \left( \frac{p}{p_1}\right) ^{\frac{\gamma-1}{2\gamma}}-1 \right],\quad if \  p\leq p_1 
		\end{cases},\\
		f_R(p)=
		\begin{cases}
			(p-p_1)\left[ \frac{A_R}{p+B_R}\right] ^{\frac{1}{2}}, \quad if \  p>p_1\\
			\frac{2a_1}{\gamma-1}\left[ \left( \frac{p}{p_1}\right) ^{\frac{\gamma-1}{2\gamma}}-1 \right],\quad if \  p\leq p_1
		\end{cases}.
	\end{aligned}
\end{equation}
$A_L,A_R,B_L,B_R$ are given by
\[
A_L=\frac{2}{(\gamma+1)\rho_1}, B_L=\frac{\gamma-1}{\gamma+1}p_1, A_R=\frac{2}{(\gamma+1)\rho_1}, B_R=\frac{\gamma-1}{\gamma+1}p_1.
\]

Experiments show that this iterative method is less sensitive to the initial value of the iteration. Lemma\ref{lemma solvability} guarantees that the solution of the iterative equation in the algorithm is unique. A classical Riemann solver $CRP(W_5,W_8)$ of the Euler equations is needed in the algorithm. The wave pattern is identical for this CRP, which consists a rarefaction wave corresponding to the $u-a$ characteristic field and a shock wave corresponding to the $u+a$ charatristic field. 

\begin{figure}[htbp]
	\label{figure test1}
	\centering
	\subfigure[density]{
		\includegraphics[width=12em]{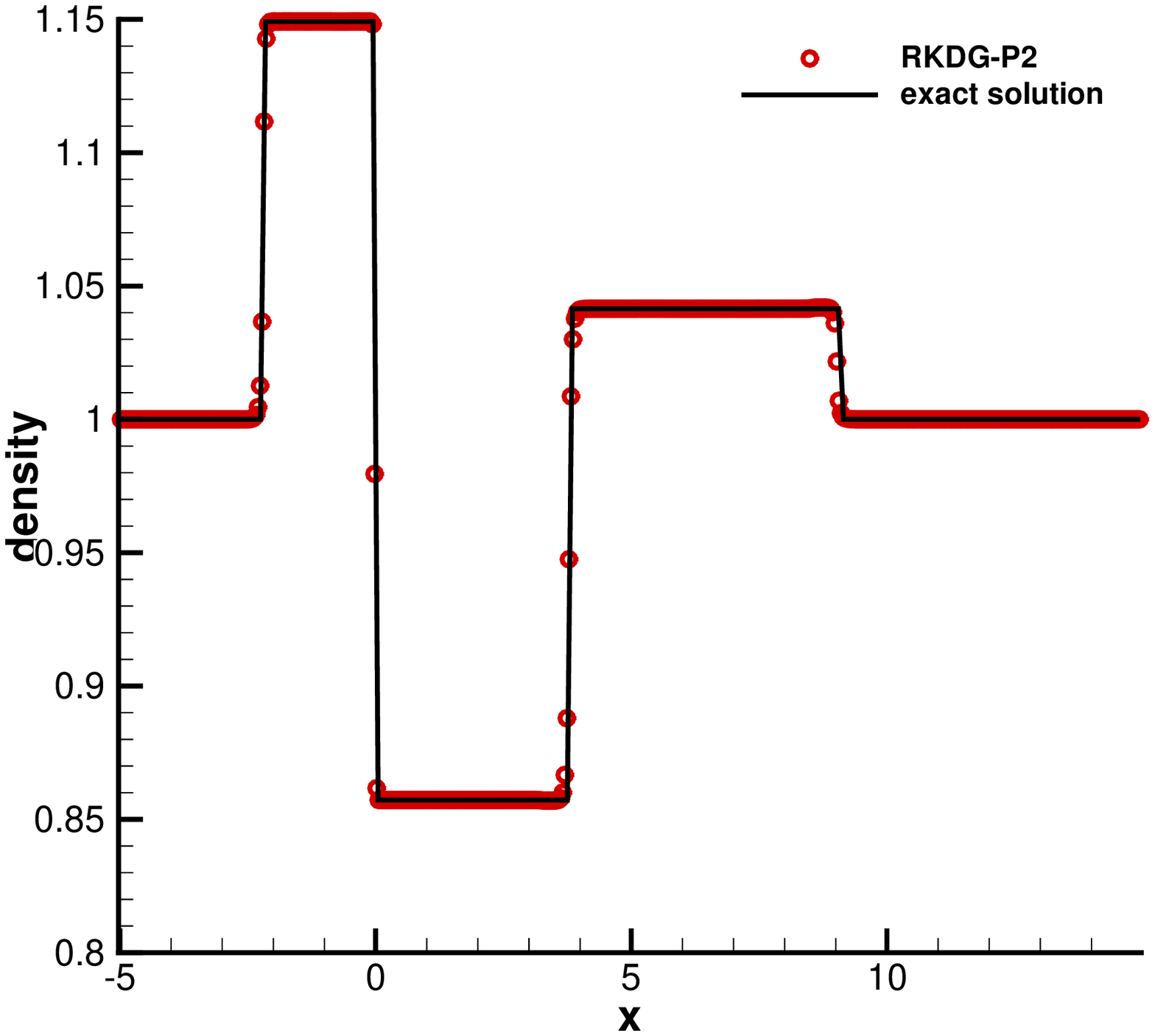}
	}
	\quad
	\subfigure[velocity]{
		\includegraphics[width=12em]{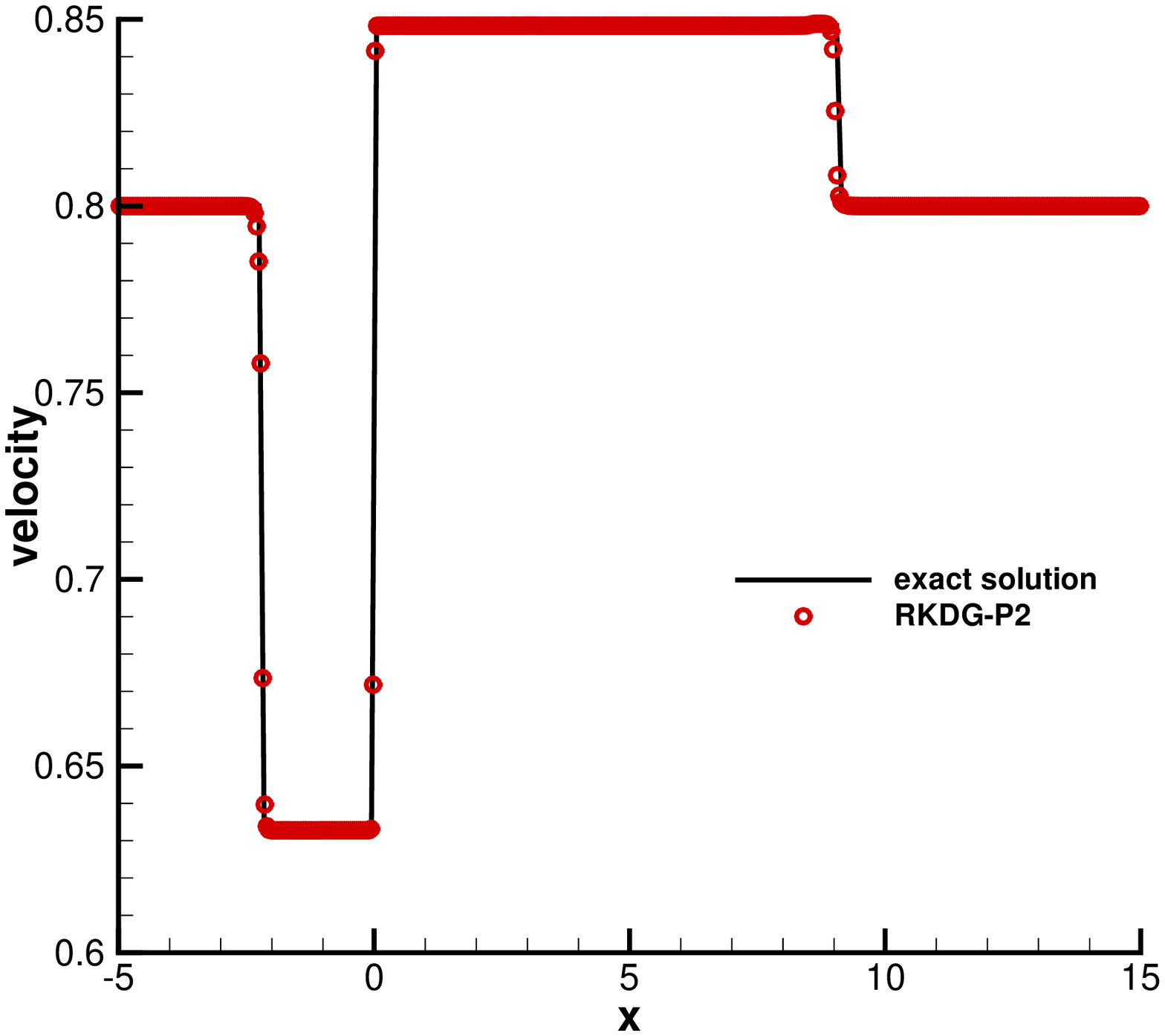}
	}
	\quad
	\subfigure[pressure]{
		\includegraphics[width=12em]{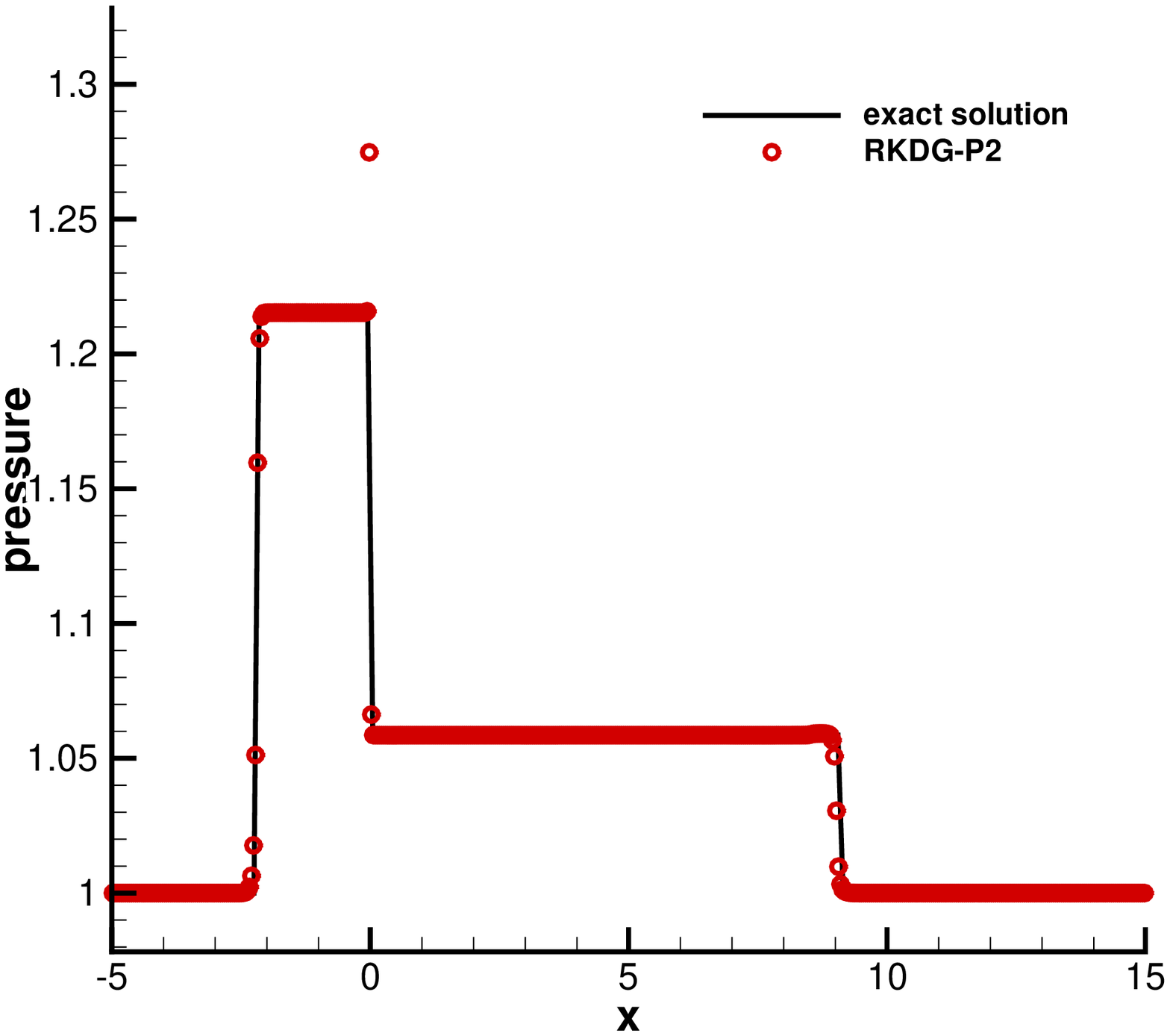}
	}
	\caption{The numerical solution obtained by the RKDG and comparion to the constructed self-similar solution for Test1 at $t=4.5s$.}
\end{figure}

\begin{figure}[htbp]
	\label{figure test2}
	\centering
	\subfigure[density]{
		\includegraphics[width=12em]{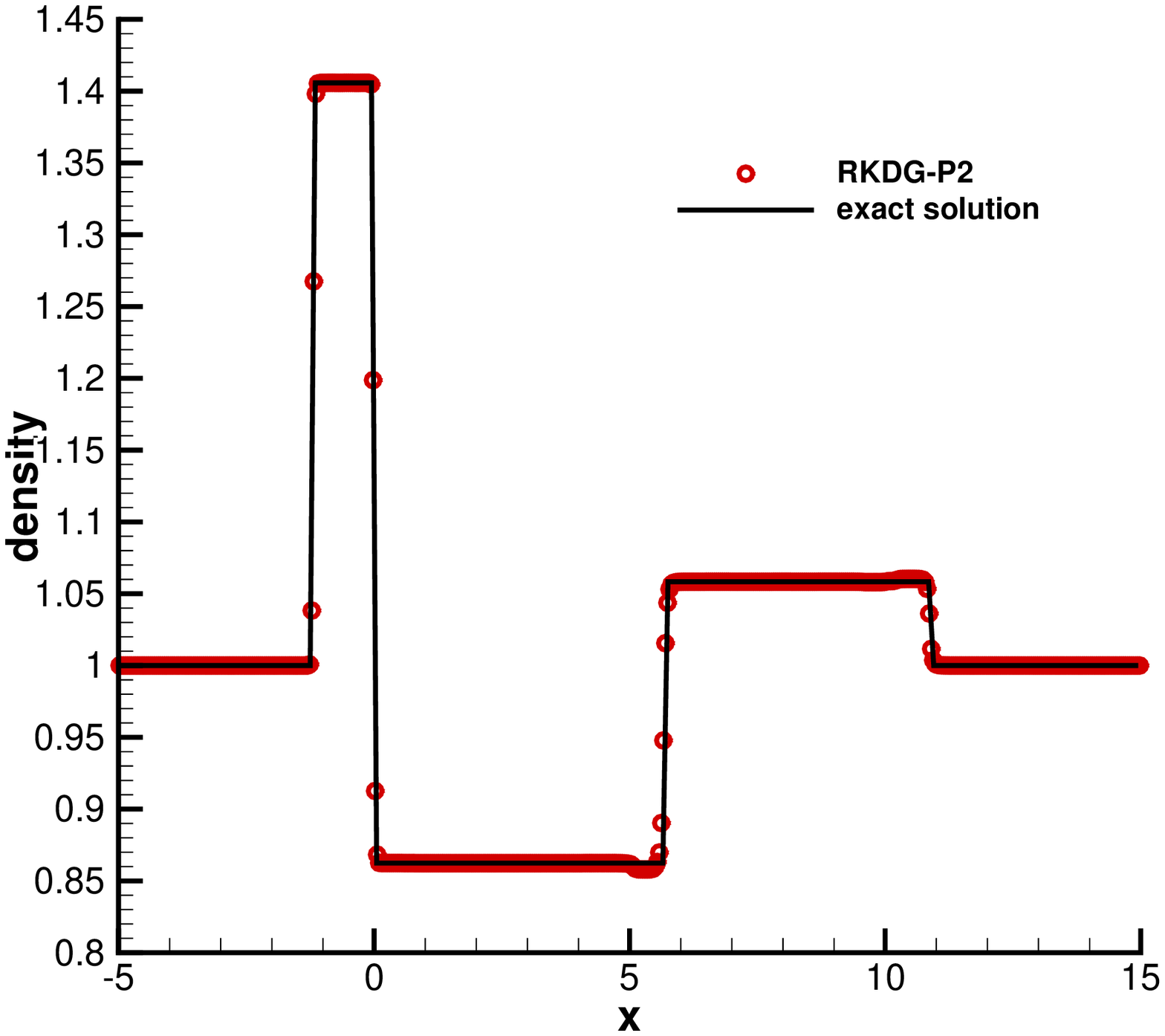}
	}
	\quad
	\subfigure[velocity]{
		\includegraphics[width=12em]{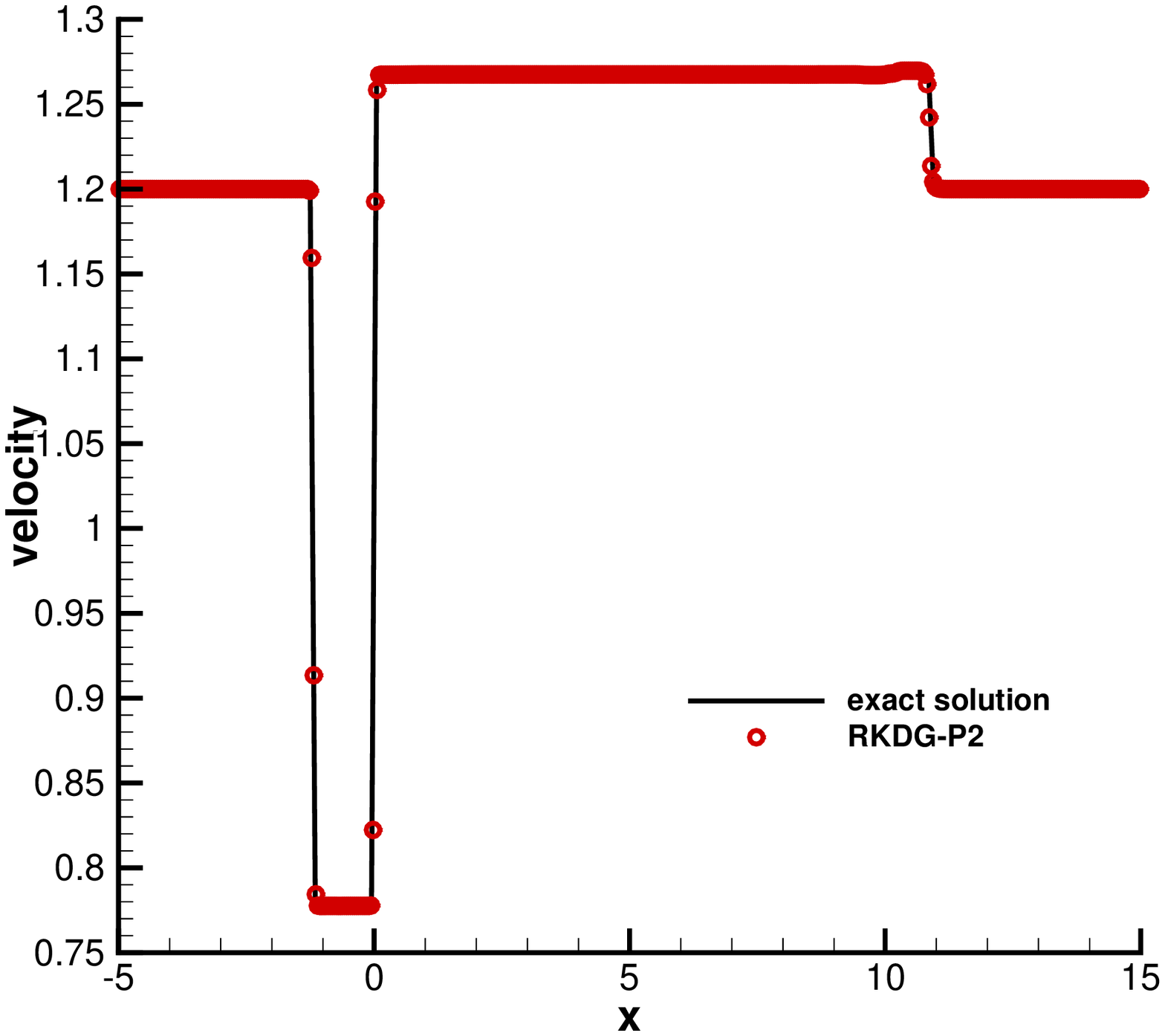}
	}
	\quad
	\subfigure[pressure]{
		\includegraphics[width=12em]{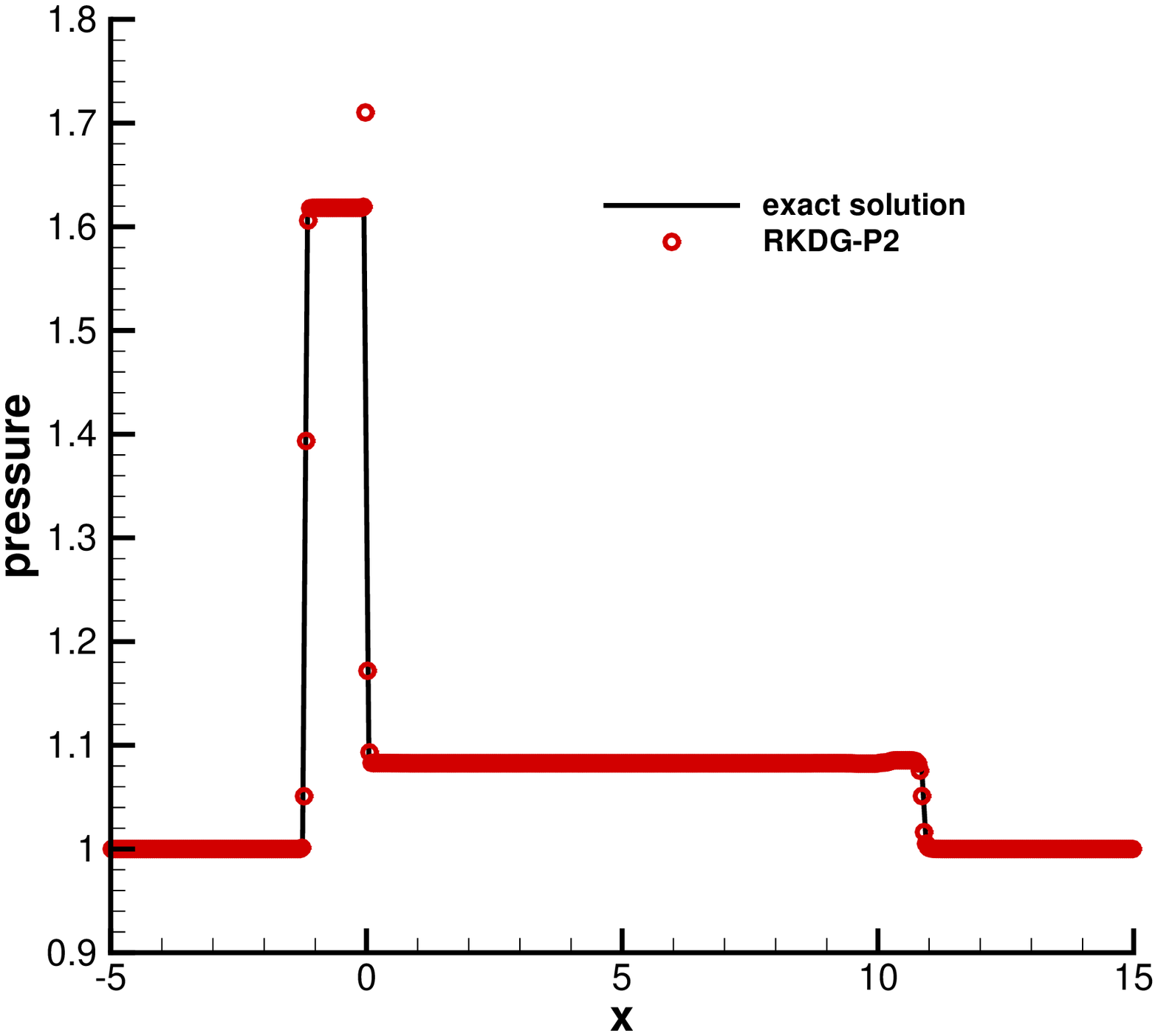}
	}
	\caption{The numerical solution obtained by the RKDG and comparion to the constructed self-similar solution for Test2 at $t=4.5s$.}
\end{figure}

\begin{figure}[htbp]
	\label{figure test3}
	\centering
	\subfigure[density]{
		\includegraphics[width=12em]{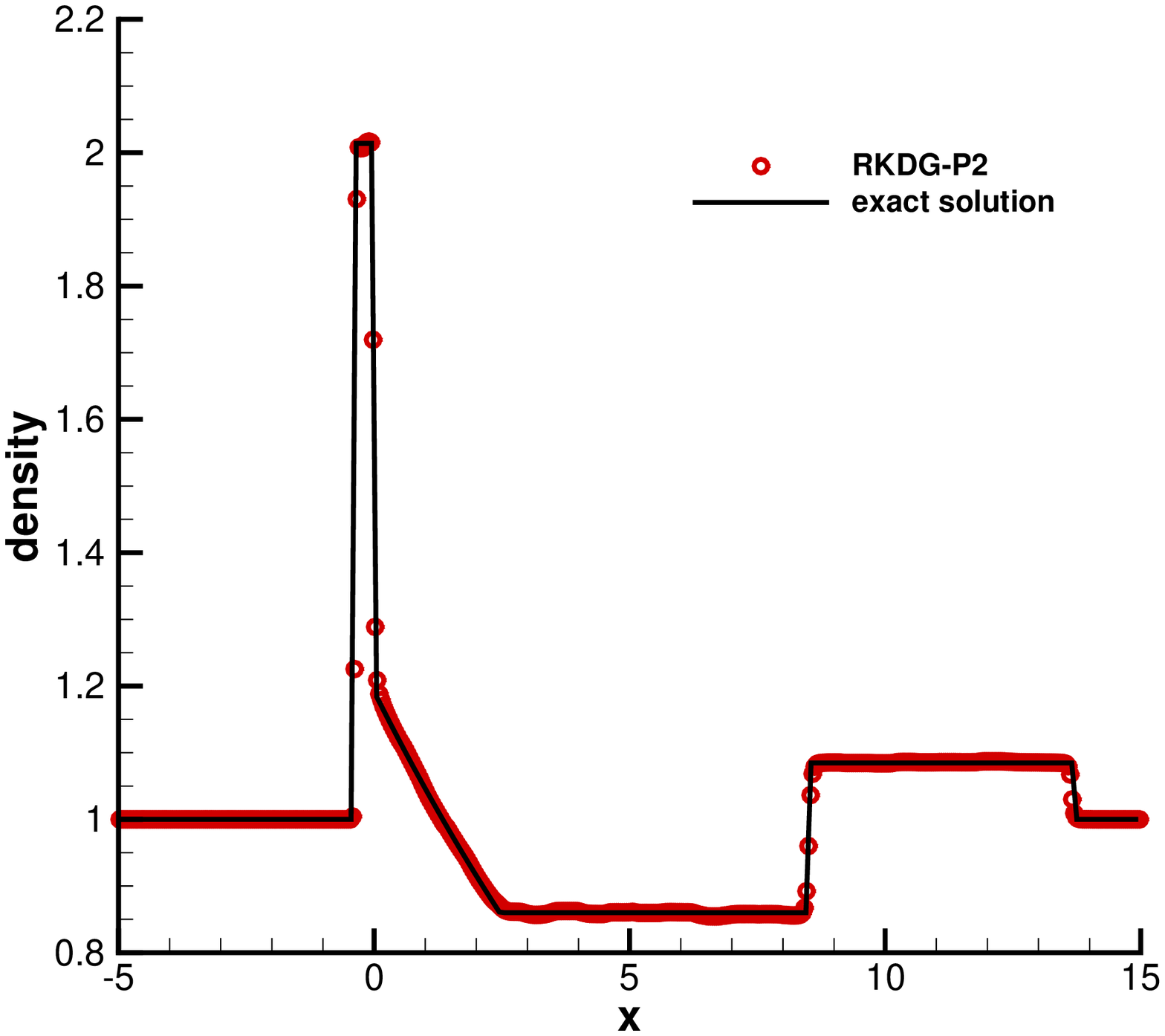}
	}
	\quad
	\subfigure[velocity]{
		\includegraphics[width=12em]{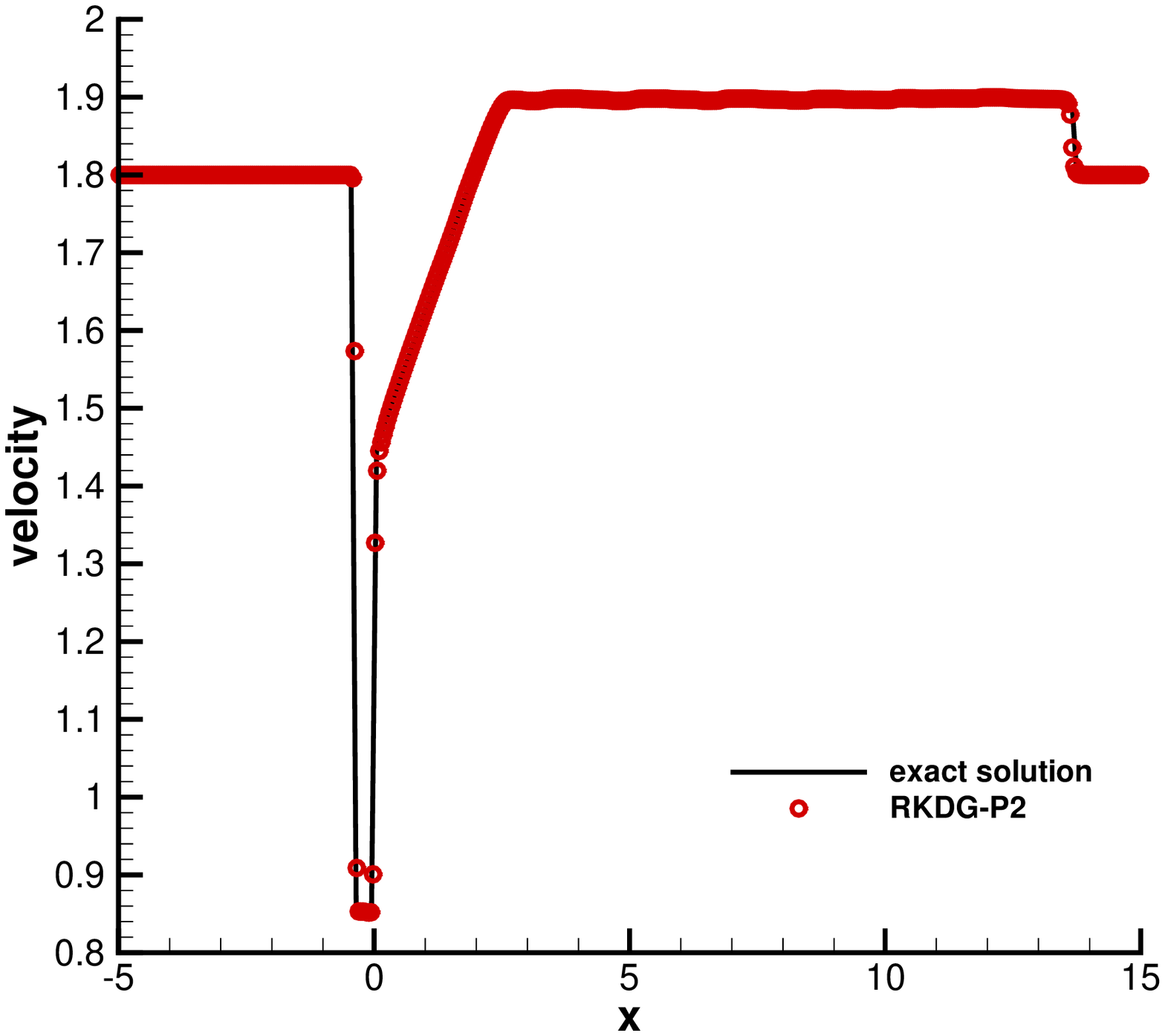}
	}
	\quad
	\subfigure[pressure]{
		\includegraphics[width=12em]{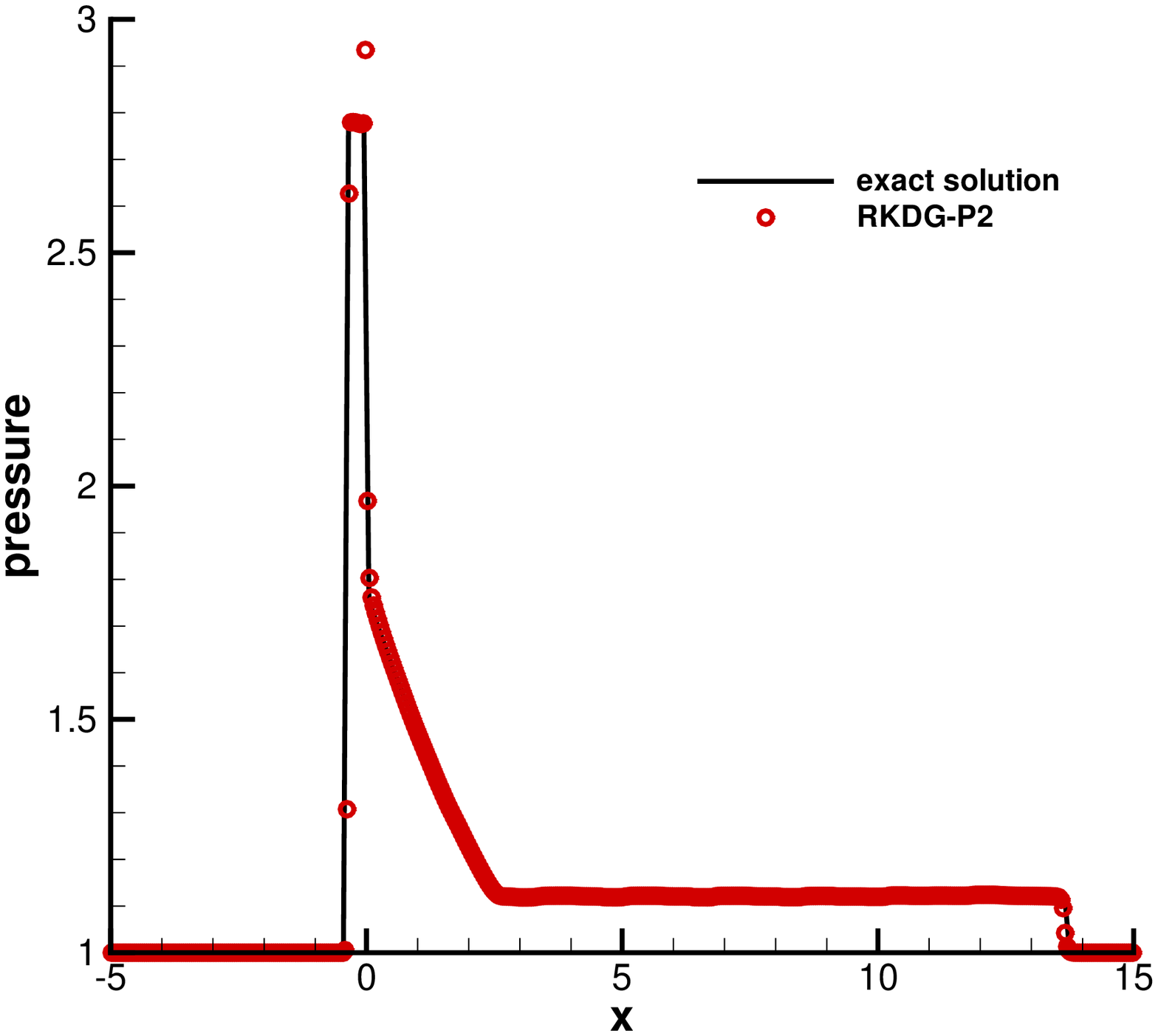}
	}
	\caption{The numerical solution obtained by the RKDG and comparion to the constructed self-similar solution for Test3 at $t=4.5s$.}
\end{figure}

\begin{figure}[htbp]
	\label{figure test4}
	\centering
	\subfigure[density]{
		\includegraphics[width=12em]{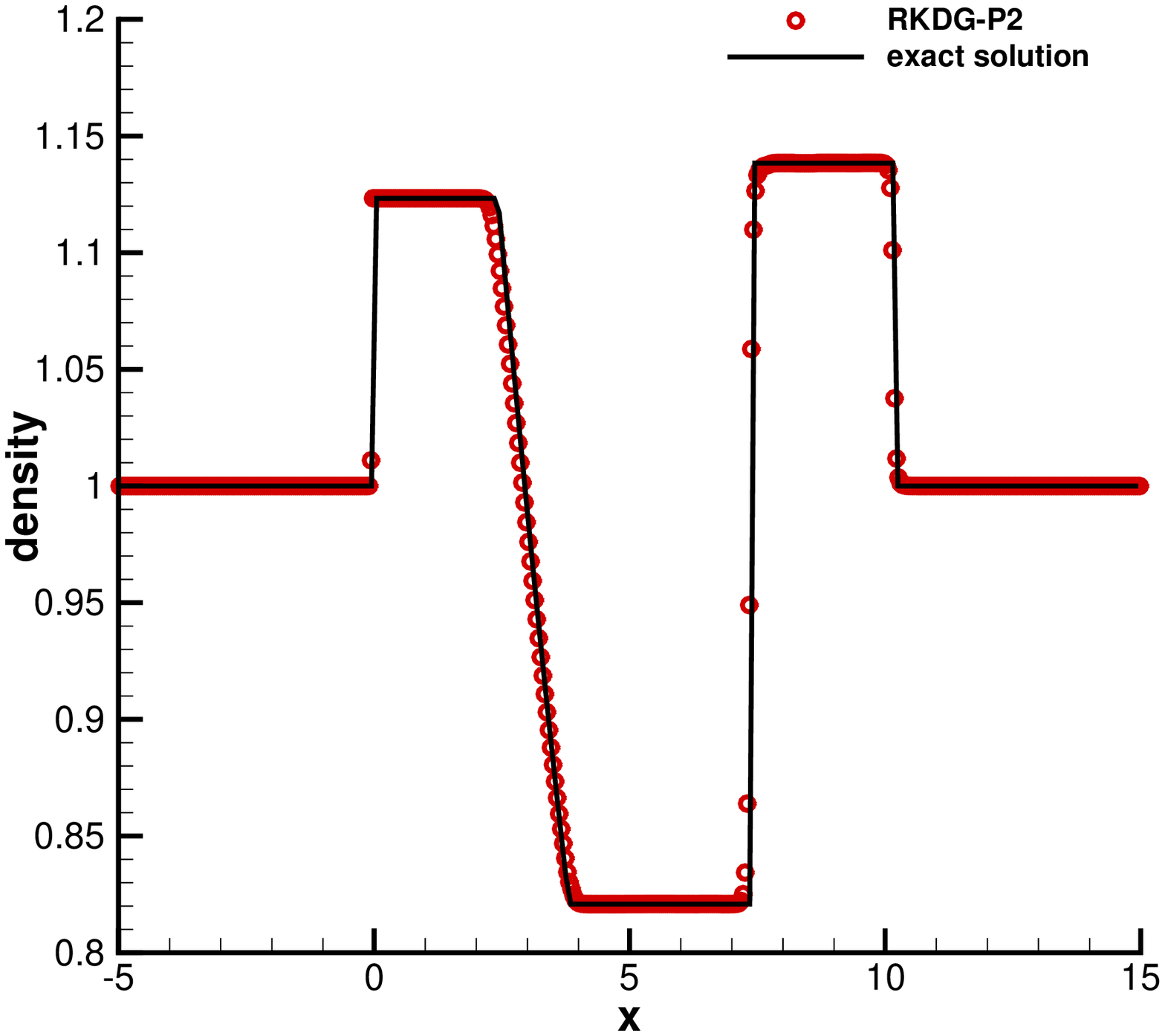}
	}
	\quad
	\subfigure[velocity]{
		\includegraphics[width=12em]{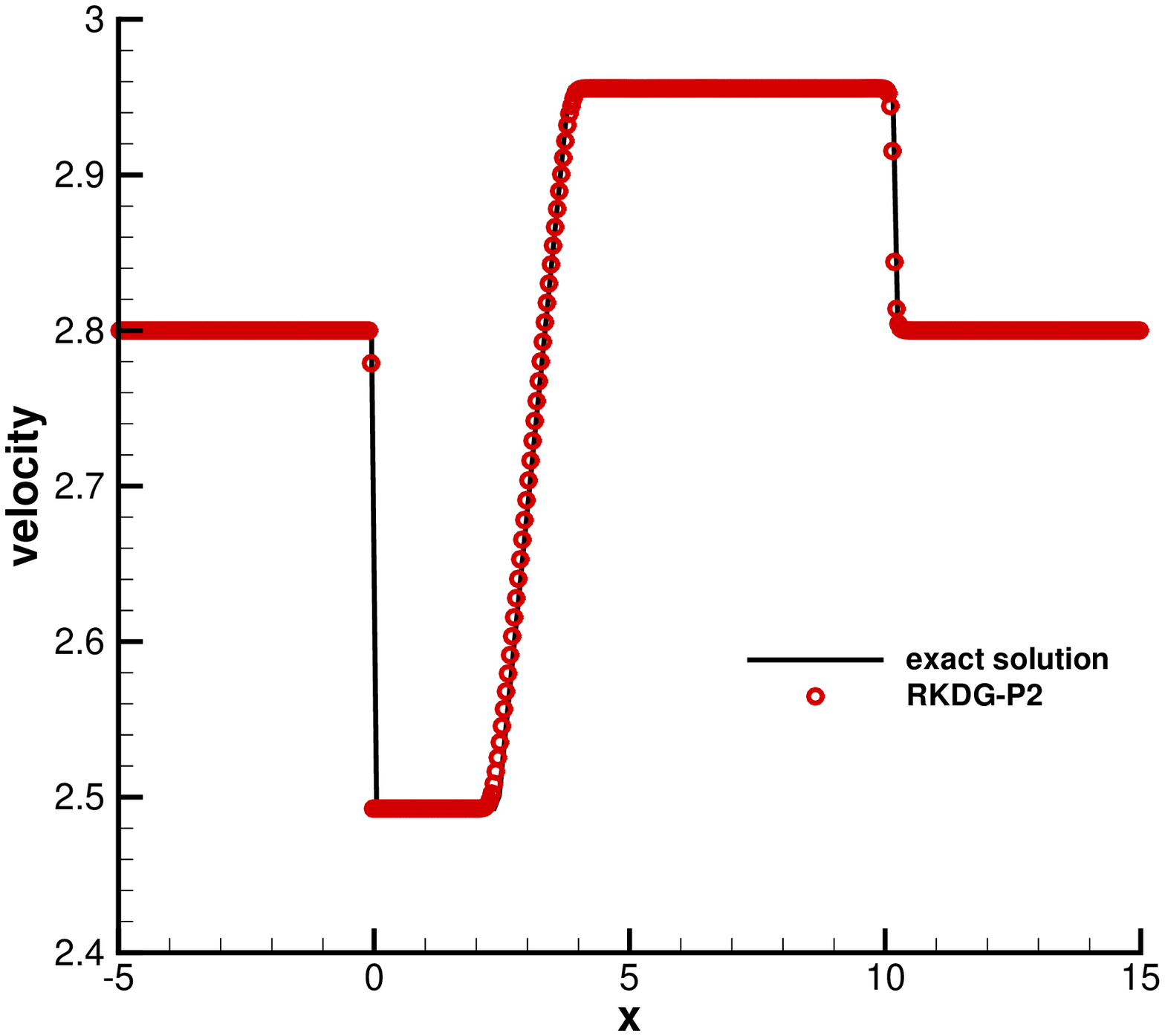}
	}
	\quad
	\subfigure[pressure]{
		\includegraphics[width=12em]{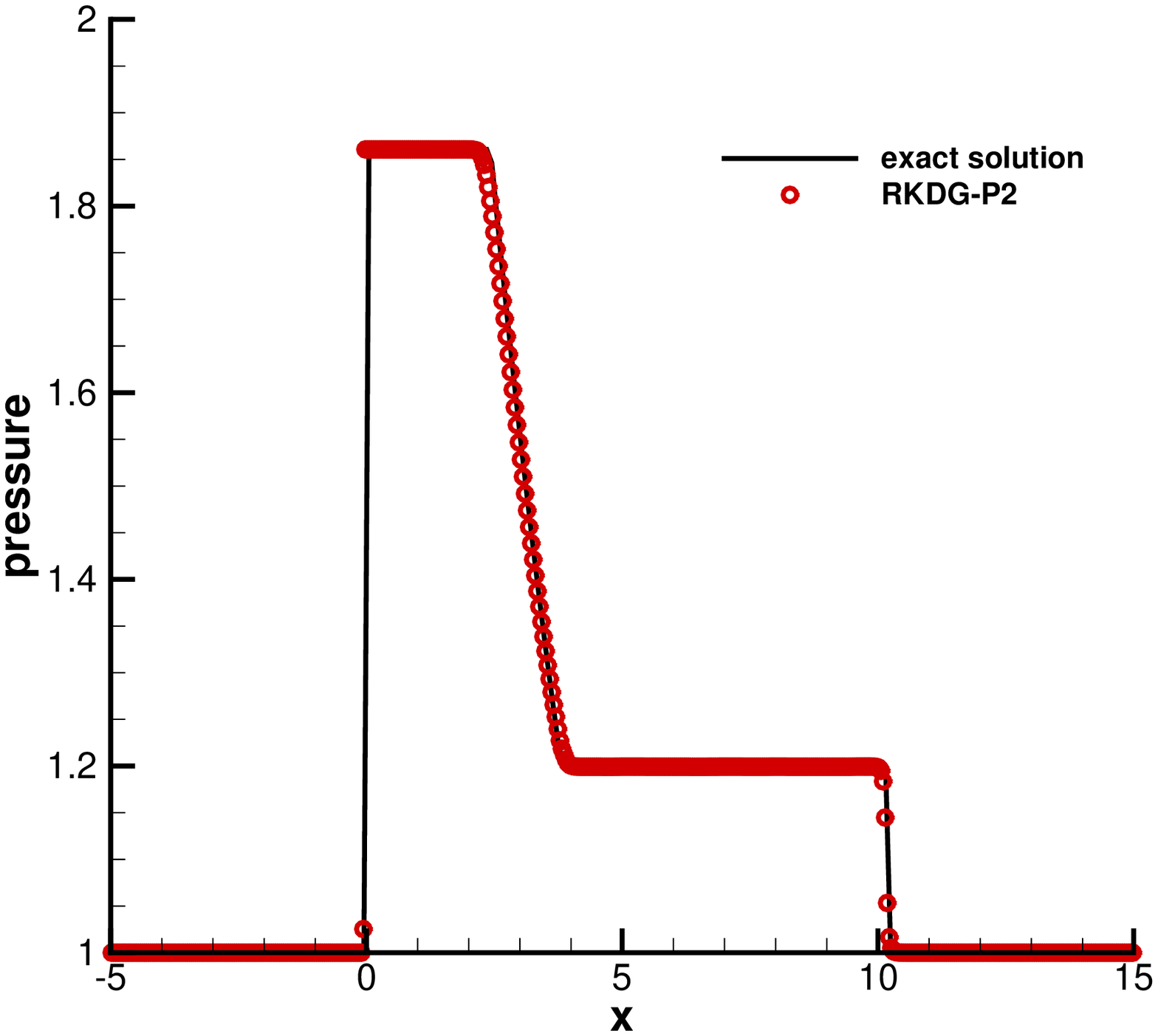}
	}
	\caption{The numerical solution obtained by the RKDG and comparion to the constructed self-similar solution for Test4 at $t=2.5s$.}
\end{figure}

\begin{figure}[htbp]
	\label{figure test5}
	\centering
	\subfigure[density]{
		\includegraphics[width=12em]{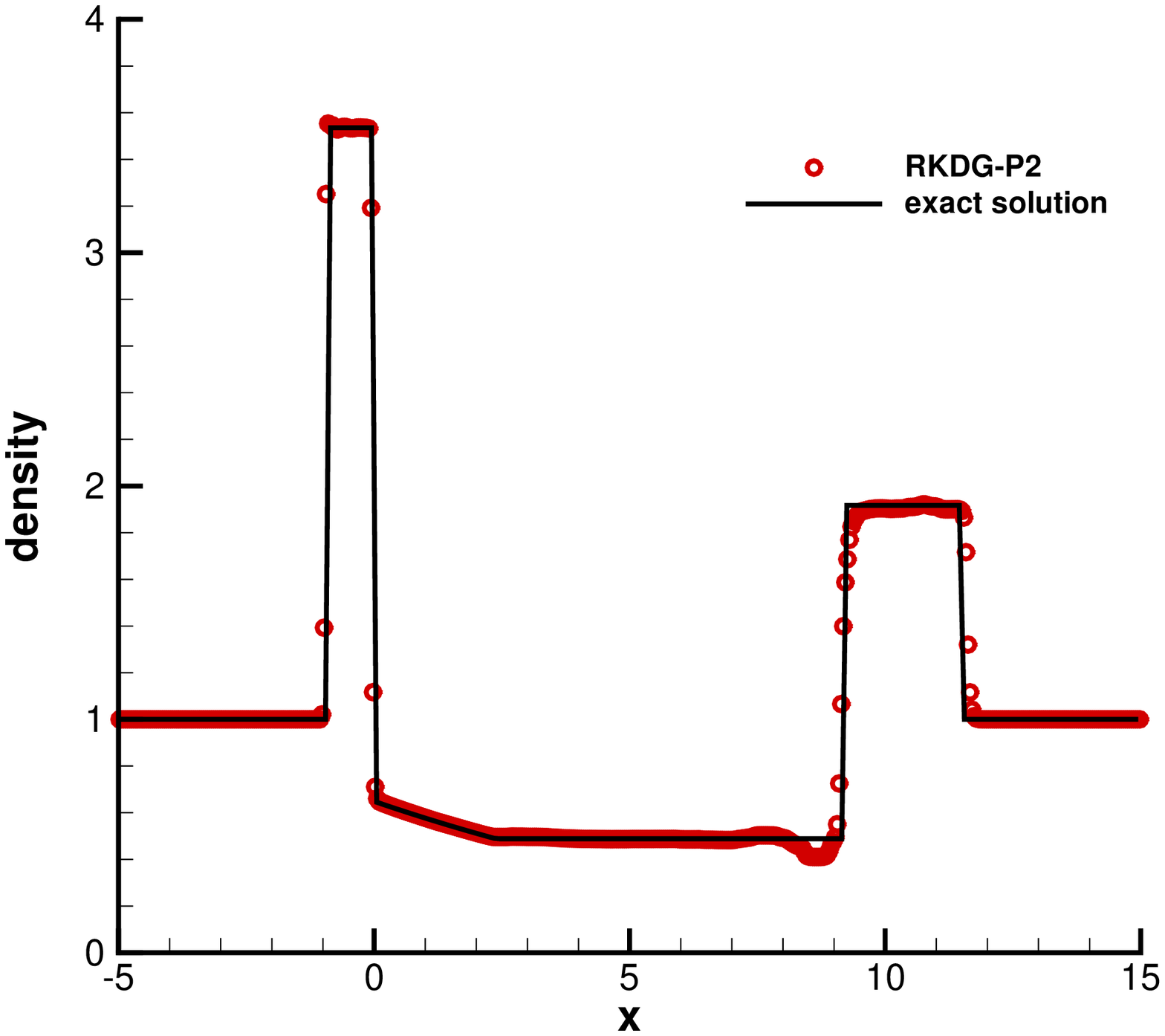}
	}
	\quad
	\subfigure[velocity]{
		\includegraphics[width=12em]{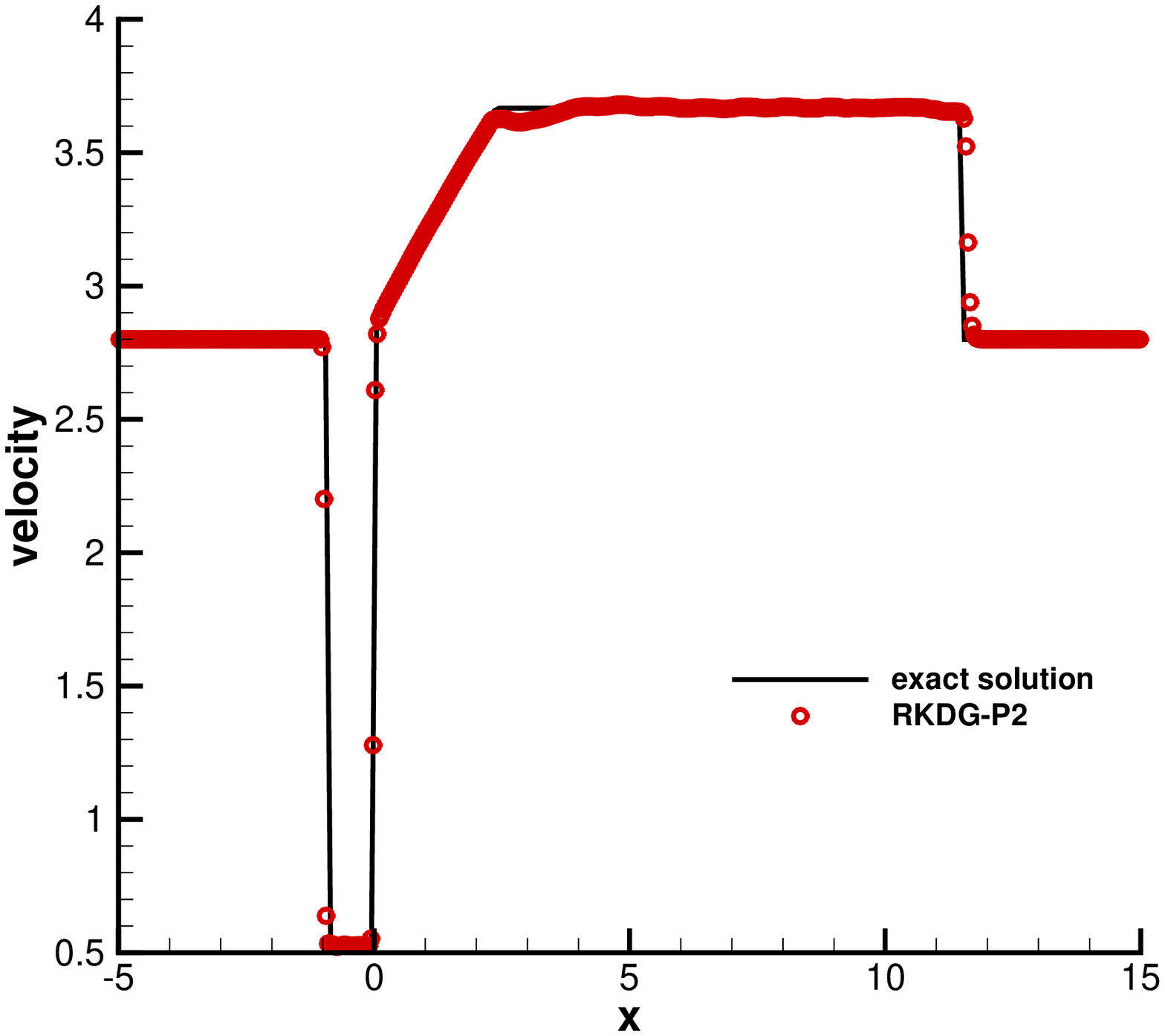}
	}
	\quad
	\subfigure[pressure]{
		\includegraphics[width=12em]{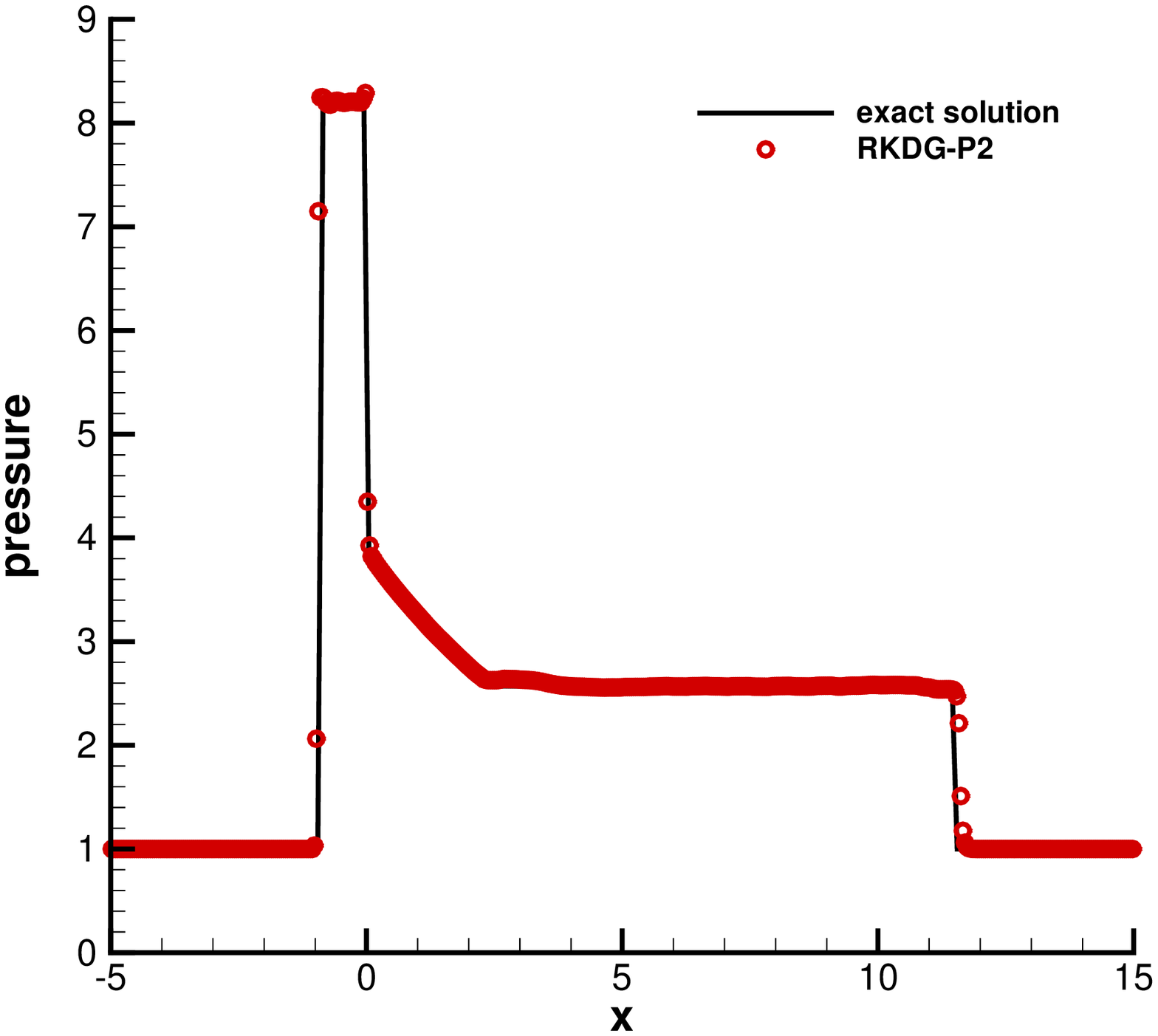}
	}
	\caption{The numerical solution obtained by the RKDG and comparion to the constructed self-similar solution for Test5 at $t=2.5s$.}
\end{figure}

We verify the existence of these three structures through numerical tests, and compare the states of the constructed self-similar solution with the numerical solution at intermediate regions. The following five tests are all for the ideal gase with $\gamma=1.4$, and their exact solutions cover three proposed structures. The initial conditions and heating parameters of the five tests are shown in Table 1.
\begin{table}[htbp]
	\label{table tests}
	\centering
	\caption{Initial conditions and parameter setting in experiments}
	\begin{tabular}{cccccc}
		\hline&
		\multicolumn{3}{c}{initial conditions}&
		\multirow{2}*{heating parameter} & 
		\multirow{2}*{solution structure} \\ \cline{2-4}
		&$\rho$&$u$&$p$& 
		\\ \hline
		Test1& 1.0&0.8&1.0&0.2&Type 1\\ \hline
		Test2& 1.0&1.2&1.0&0.2&Type 1\\ \hline
		Test3& 1.0&1.8&1.0&0.2&Type 2\\ \hline
		Test4& 1.0&2.8&1.0&0.2&Type 3\\ \hline
		Test5& 1.0&2.8&1.0&2.0&Type 2\\ \hline
	\end{tabular} 
\end{table}
The comparison between the numerical solutions and the exact solutions are shown in Figure\ref{figure test1,figure test2,figure test3,figure test4,figure test5}. In the first four tests, the heating parameter is $k=0.2$, and $k(\gamma^2-1)<1$ holds. $M_{*}$ and $M_{**}$ are 0.6136 and 1.8130, respectively. The root of $Y(M_1)$ is 1.0620, which is less than $M_{**}$, hence the solution is unique. In the last test, the value of $k$ is 2.0. $k(\gamma^2-1)>1$ holds at this time. According to Theorem\ref{theorem uniqueness}, the exact solutions of these five Riemann problem are all unique.

We apply the Runge-Kutta discontinuous Galerkin method (RKDG) for numerical simulation (see \cite{yang2013discontinuous}). In the spatial direction, the solution is discretized by piece-wise second-order ($P2$) polynomials. A total variation diminishing (TVD) limiter is empolyed to avoid numerical oscillations. In the time direction, we use the third-order Runge-Kutta method. The source is processed through splitting. Note that the source term should be precessed at each time step of Runge-Kutta method.

It can be observed from figures that the constructed exact solution and numerical solution fit well in each intermediate regions. While inside several cells near origin, the solutions of pressure are some different between the exact solution and the numerical solution. The reason is that directly splitting the source term has no "well-balanced" property.

\section{Conclusions}
\label{sec:conclusions}

This paper focused on the Riemann problem of the Euler equations with a Dirac delta-source in the energy conservation eqution. The double CRPs frame was proposed to construct the self-similar solutions. We proved that there are three types of the solution and studied the uniqueness of the solutions under this frame. The present frame is completely different from the existing method of dealing with the Riemann problem with discontinuous source. Based on the solution of CRP, the strategy of the present frame is the elimination of unconscionable waves on the general structure (Figure\ref{figure general structure}). This frame can be applied to the Riemann problem of other hyperbolic systems with Dirac delta-function sources or other sources. Compared with some existing methods, it is more simple and can naturally cover all possible structures. We verified the double CRPs frame by comparing the constructed self-similar solution with the numerical solution obtained by RKDG. Besides, the constructed solutions can be used to evaluate the existing numerical methods for source terms, such as \cite{greenberg1997analysis,jin2005two,kroner2005numerical}.

The uniqueness of the self-similar solutions under the double CRPs frame for arbitrary initial conditions is an open question, which is the goal of our future work. In addition, future work should focous on the heating addition of unsteady flow, in which the transition of the three structures proposed in this paper may occur.

\section*{Appendix}
\subsection*{Proof of Lemma\ref{lemma uniqueness1}}
\label{proof uniqueness1}
\mbox{}\\
\begin{Proof}
	For the solution of Type 1, (\ref{Lshock formula}) and (\ref{heat addition formula}) imply
	\[
	\begin{aligned}
		&p_4/p_1=f_2(M_1,M_{SL}),\\
		&p_8/p_7=p_1/p_4\times p_4/p_5=\left[ f_2(M_1,M_{SL}\psi(f_5(M_1,M_{SL})))\right] ^{-1},
	\end{aligned}
	\]
	where $M_{SL}$ is obtained from equation (\ref{left equations}) and is determined by $M_1$.
	
	For Structure 2, we have
	\[
	p_4/p_1=f_2(M_1,M_{SL}),
	\]
	where $M_{SL}$ is obtained from equation (\ref{structure2 mach number}) and is determined by $M_1$.
	
	From the elementary wave equation in \cite{0Riemann}, we have
	\[
	\begin{aligned}
		&u_6=u_5-\frac{2a_5}{\gamma-1}\left[ \left( \frac{p_6}{p_5}\right) ^{\frac{\gamma-1}{2\gamma}}-1\right], \\
		&u_7=u_8+\sqrt{2}(p_7-p_8)\left[ (\gamma+1)p_7\rho_8+(\gamma-1)p_8\rho_8\right] ^{-\frac{1}{2}}.
	\end{aligned}
	\]
	According to $u_6=u_7$, $p_6=p_7$ and $U_1=U_8$, it follows that
	\[
	u_5-\frac{2a_5}{\gamma-1}\left[ \left( \frac{p_6}{p_5}\right) ^{\frac{\gamma-1}{2\gamma}}-1\right]
	=u_1+\sqrt{2}(p_6-p_1)\left[ (\gamma+1)p_6\rho_1+(\gamma-1)p_1\rho_1\right] ^{-\frac{1}{2}}.
	\]
	Divide both sides of the above equation by $a_1$, then
	\[
	\frac{u_5}{u_1}M_1-\frac{2}{\gamma-1}\frac{M_1}{M_5}\frac{u_5}{u_1}\left[ \left( \frac{p_6}{p_5}\right) ^{\frac{\gamma-1}{2\gamma}}-1\right]
	=M_1+\frac{\sqrt{2}}{a_1}(p_6-p_1)\left[ (\gamma+1)p_6\rho_1+(\gamma-1)p_1\rho_1\right] ^{-\frac{1}{2}}.
	\]
	Apply the equation of state to the above equation, we have
	\[
	\frac{u_5}{u_1}M_1-\frac{2}{\gamma-1}\frac{M_1}{M_5}\frac{u_5}{u_1}\left[ \left( \frac{p_6}{p_5}\right) ^{\frac{\gamma-1}{2\gamma}}-1\right]
	=M_1+\sqrt{\frac{2}{\gamma}}(\frac{p_6}{p_1}-1)\left[ (\gamma+1)\frac{p_6}{p_1}+(\gamma-1)\right] ^{-\frac{1}{2}}.
	\]
	Therefore
	\[
	\left( A\left( \frac{p_6}{p_5}\right) ^{\frac{\gamma-1}{2\gamma}}+B\right)\sqrt{C\frac{p_6}{p_5}+\gamma-1}
	-\sqrt{\frac{2}{\gamma}}\left( \frac{C}{\gamma+1}\frac{p_6}{p_5}-1\right)=0 ,
	\]
	where $A=-\frac{2}{\gamma-1}M_1\phi(M_*)f_4(M_1,M_{SL})$, $B=\frac{\gamma+1}{\gamma-1}M_1\phi(M_*)f_4(M_1,M_{SL})$, \\	
	$C=(\gamma+1)\psi(M_*)f_2(M_1,M_{SL})$ and $M_{SL}$ is obtained from (\ref{structure2 mach number}). Thus the strength of $WR_1$ in the solution of Type 2 is determined by $M_1$.
	
	For $WR_3$, we have
	\[
	\frac{p_8}{p_7}=\left[ f_2(M_1,M_{SL})\psi(M_*)\frac{p_6}{p_5}\right] ^{-1}.
	\]
	Therefore the strength of each nonlinear wave in the solution of Type 2 is determined by $M_1$.
	
	The analyse of Type 3 is similar, as follows.
	\[
	\begin{aligned}
		&\left( A'\left( \frac{p_6}{p_5}\right) ^{\frac{\gamma-1}{2\gamma}}+B'\right)\sqrt{C'\frac{p_6}{p_5}+\gamma-1}
		-\sqrt{\frac{2}{\gamma}}\left( \frac{C'}{\gamma+1}\frac{p_6}{p_5}-1\right)=0, \\
		&\frac{p_8}{p_7}=\left[ \frac{p_5}{p_1}\frac{p_6}{p_5}\right] ^{-1},
	\end{aligned}
	\]
	where $A'=-\frac{2}{\gamma-1}\frac{M_1}{M_5}\frac{u_5}{u_1}$, $B'=\frac{u_5}{u_1}M_1+\frac{2}{\gamma-1}\frac{M_1}{M_5}\frac{u_5}{u_1}$, $C'=(\gamma+1)\frac{p_5}{p_1}$. $\frac{u_5}{u_1}$, $\frac{p_5}{p_1}$ and $M_5$ in the above equation are obtained from Lemma\ref{lemma proof formula1} and they are all detrmined by $M_1$.
\end{Proof}

\subsection*{Proof of Lemma\ref{lemma uniqueness2}}
\label{proof uniqueness2}
\begin{Proof}
	From (\ref{left equation}), we have 
	\[
	X(M_1,M_4)=0.
	\]
	Some tedious manipulation yields
	\[ \frac{\partial f_5}{\partial M_{SL}}>0, \]
	and
	\[ \frac{\partial f_6}{\partial M_{SL}}=\left( \frac{\partial f_5}{\partial M_{SL}}\right) ^{-1}>0. \]
	
	Then
	\begin{equation}
		\label{uniqueness structure1 formula}
		\frac{\partial X}{\partial M_4}=
		M_1\frac{\partial}{\partial M_4}(f_4(M_1,M_{SL})\phi(M_4))
		-\sqrt{\frac{\beta}{\gamma}}\frac{\partial }{\partial M_4}
		\left(\frac{f_2(M_1,M_{SL})\psi(M_4)-1}{\sqrt{1+\tau f_2(M_1,M_{SL})\psi(M_4)}}\right).
	\end{equation}
	
	For the first part on the right of (\ref{uniqueness structure1 formula}), we have
	\[
	\frac{\partial}{\partial M_4}(f_4(M_1,M_{SL})\phi(M_4))=\frac{\partial f_4(M_1,M_4)}{\partial M_4}\phi(M_4)+\frac{d \phi(M_4)}{d M_4}f_4(M_1,M_{SL}),
	\]
	\[
	\frac{\partial f_4(M_1,M_4)}{\partial M_{4}}=
	\frac{\partial f_4(M_1,M_4)}{\partial M_{SL}}\frac{\partial M_{SL}}{\partial M_{4}}=
	\frac{2(1+(M_1-M_{SL})^{-2})}{(\gamma+1)M_1}\frac{\partial f_6}{\partial M_4}>0,
	\]
	\[
	\phi ' (M_4)>0.
	\]
	Thus
	\[
	\frac{\partial}{\partial M_4}(f_4(M_1,M_{SL})\phi(M_4))>0.
	\]
	
	For the second part on the right of (\ref{uniqueness structure1 formula}), we have
	\[
	\frac{\partial }{\partial M_4}
	\left(\frac{f_2(M_1,M_{SL})\psi(M_4)-1}{\sqrt{1+\tau f_2(M_1,M_{SL})\psi(M_4)}}\right)
	=\frac{\tau (f_2(M_1,M_{SL})\psi(M_4)+1)+2}{2(f_2(M_1,M_{SL})\psi(M_4))^{3/2}}
	\frac{\partial}{\partial M_4}(f_2(M_1,M_{SL})\psi(M_4)),
	\]
	\[
	\frac{\partial}{\partial M_4}(f_2(M_1,M_{SL})\psi(M_4))
	=-\left( \frac{\partial f_2}{\partial M_1}-\frac{\partial f_2}{\partial M_{SL}}\right) \psi(M_4)+\frac{d \psi(M_4)}{d M_{SL}}f_2(M_1,M_{SL}),
	\]
	\[
	\frac{\partial f_2}{\partial M_1}-\frac{\partial f_2}{\partial M_{SL}}>0,
	\]
	\[
	\psi ' (M_4)<0.
	\]
	Thus
	\[
	\frac{\partial}{\partial M_4}(f_2(M_1,M_{SL})\psi(M_4))<0,
	\]
	and
	\[
	\frac{\partial }{\partial M_4}
	\left(\frac{f_2(M_1,M_{SL})g_2(M_1,M_{SL})-1}{\sqrt{1+\tau f_2(M_1,M_{SL})g_2(M_1,M_{SL})}}\right)<0.
	\]
	
	Note that we have actually proved that
	\[
	\frac{\partial X}{\partial M_4}>0.
	\]
	
	According to Lemma\ref{lemma upwind condition}, it follows that $M_4<M_*$. Consequently
	\[Y(M_1)=X(M_1,M_*)\geq X(M_1,M_4)=0. \]
\end{Proof}

\begin{Proof}
	The solution right to the t-axis of Type 2 is the solution of $CRP(U_5,U_8)$, whose wave pattern is $R(U_5,U_6)\oplus C(U_6,U_7)\oplus S(U_7,U_8)$, it follows that (see \cite{liu2005the})
	\[u_5-u_8\leq \sqrt{\frac{\beta p_8}{\rho_8}}\frac{p_5/p_8-1}{\sqrt{1+\tau p_5/p_8}}. \]
	Represent the above equation by the Mach numbers, as follows.
	\[ M_1(f_4(M_1,M_{SL})\phi(M_*)-1)
	-\sqrt{\frac{\beta}{\gamma}}\frac{f_2(M_1,M_{SL})\psi(M_*)-1}{\sqrt{1+\tau f_2(M_1,M_{SL})\psi(M_*)}}\leq 0. \]
	Thus
	\[
	Y(M_1)\leq 0.
	\]
	We denote
	\[
	M_{SL}'=\frac{s_L}{a_4},
	\]
	where $s_L$ is the speed of $WL_1$ and $a_4$ is the speed of sound in region 4.	Similar to (\ref{Lshock formula}), it holds that
	\[
	M_1=\frac{((\gamma-1)M_4+2M_{SL}')(M_4-M_{SL}')+2}{\sqrt{2\gamma(\gamma-1)(M_4-M_{SL}')^4
			+(6\gamma-\gamma^2-1)(M_4-M_{SL}')^2-2(\gamma-1)}}.
	\]
	Define
	\[
	g(M_4,M_{SL}')=\frac{((\gamma-1)M_4+2M_{SL}')(M_4-M_{SL}')+2}{\sqrt{2\gamma (M_4-M_{SL}')^2-\gamma+1}\sqrt{(\gamma-1)(M_4-M_{SL}')^2+2}}.
	\]
	The domain of $g$ satisfies
	\[
	M_4\leq M_*,\quad M_{SL}'\leq 0.
	\]
	A routine computation gives rise to
	\[
	\frac{\partial g}{\partial M_{SL}'} >0.
	\]
	Then
	\[
	g(M_4,M_{SL}')\leq g(M_4,0)=\sqrt{\frac{(\gamma-1)M_4^2+2}{2\gamma M_4^2-\gamma+1}},
	\]
	and
	\[
	M_1=g(M_*,M_{SL}')\leq g(M_*,0)=\sqrt{\frac{(\gamma-1)M_*^2+2}{2\gamma M_*^2-\gamma+1}}.
	\]
	According to Lemma\ref{lemma prandtl relation}, it follows that
	\[ M_1\leq M_{**}. \]
	
	The proof of (iii) is obvious.
\end{Proof}

\section*{Acknowledgments}
This work was supported by the NSFC-NSAF joint fund [No. U1730118]; and the Science Challenge Project [No. JCKY2016212A502].

\bibliographystyle{plain}
\bibliography{references}

\end{document}